\documentclass[12pt,twoside,leqno]{article}
\usepackage{amsmath}
\usepackage{amssymb}
\usepackage{amsxtra}
\usepackage{amscd}
\usepackage{amsthm}
\usepackage[mathscr]{eucal}
\usepackage{xr}
\usepackage{color}

\swapnumbers

\theoremstyle{plain}
\newtheorem{thm}[subsection]{Theorem}
\newtheorem{prop}[subsection]{Proposition}

\newtheorem{lem}[subsection]{Lemma}

\theoremstyle{definition}
\newtheorem{defn}[subsection]{Definition}

\newtheorem{rem}[subsection]{Remark}
\newtheorem{para}[subsection]{}

\newtheorem{sbrem}[subsubsection]{Remark}

\newenvironment{pf}{\proof[\proofname]}{\endproof}
\newenvironment{pf*}[1]{\proof[#1]}{\endproof}

\newcommand{\N}{{\mathbb{N}}}
\newcommand{\Q}{{\mathbb{Q}}}
\newcommand{\Z}{{\mathbb{Z}}}
\newcommand{\R}{{\mathbb{R}}}

\newcommand{\isom}{\cong}
\newcommand{\Spec}{\operatorname{Spec}}

\newcommand{\fs}{{\rm fs \ }}
\newcommand{\gp}{\mathrm{gp}}
\newcommand{\Gm}{\mathbb{G}_m} %
\newcommand{\Ga}{\mathbb{G}_a} %
\newcommand{\Gmlog}{\mathbb{G}_{m,\log}}

\renewcommand{\O}{{\mathcal{O}}}

\newcommand{\id}{\mathrm{id}}
\renewcommand{\tilde}{\widetilde}
\newcommand{\Hom}{\operatorname{Hom}}
\newcommand{\End}{\operatorname{End}}
\newcommand{\Ext}{\operatorname{Ext}}

\newcommand{\Biext}{\operatorname{Biext}}
\newcommand{\BIEXT}{\operatorname{\mathcal{B}\it{iext}}}

\newcommand{\HOM}{\operatorname{\mathcal{H}\it{om}}}
\newcommand{\Aut}{\operatorname{Aut}}

\newcommand{\Ker}{\operatorname{Ker}}

\renewcommand{\bar}[1]{\overline{#1}}

\def\spcheck{^\vee}

\def \llet {\mathrm {l\acute{e}t}} %
\def \et {\mathrm {\acute{e}t}}

\newcommand{\cl}{Claim}
\theoremstyle{definition}

\newtheorem*{clm*}{\cl}

\theoremstyle{plain}

\def \fs {\mathrm {fs}}
\def \ket {\mathrm {k\acute{e}t}}

\def \overc#1{\overset {\lower 0.3ex \hbox{${\;}_{\circ}$}}{#1}}

\let\refsave=\ref
\def\ref#1{\textup{\refsave{#1}}}

\newcommand{\sym}{\mathrm{sym}}

\newcommand{\oc}[1]{{\overc#1}}

\newcommand{\upc}{\overset{\circ}\to}
\newcommand\Cal{\mathcal}
\newcommand\define{\newcommand}
\renewcommand\bold{\Bbb}

\define\bZ{\bold Z}

\define\bR{\bold R}
\define\bQ{\bold Q}
\define\bN{\bold N}
\define{\cS}{\Cal S}
\define{\Lie}{\mathrm{Lie}\,} %
\define{\coLie}{\mathrm{coLie}\,} %
\define{\cH}{\Cal H}
\define{\cExt}{{\Cal E}xt}
\define{\cHom}{{\Cal H}om}
\define{\cO}{\Cal O}
\define{\an}{\mathrm{an}}

\newcommand{\ol}[1]{\overline{#1}}

\def \upcf {\overset {\lower 0.3ex \hbox{${\;}_{\circ}$}} f}
\def \upcp {\overset {\lower 0.3ex \hbox{${\;}_{\circ}$}} p}
\def \upc#1{\overset {\lower 0.3ex \hbox{${\;}_{\circ}$}}{#1}}

\newcommand{\Sig}{\Sigma}
\newcommand{\sig}{\sigma}

\newcommand{\cA}{\mathcal A}
\newcommand{\cB}{\mathcal B}
\newcommand{\cC}{\mathcal{C}}

\newcommand{\compo}{{\circ}}
\newcommand{\cI}{\mathcal I}
\newcommand{\cT}{\mathcal T}

\begin{document}
\title
{Logarithmic abelian varieties, \\ Part VII: Moduli} 
\author{Takeshi Kajiwara, Kazuya Kato, and Chikara Nakayama}
\date{}
\maketitle
\centerline{\it{Dedicated to Professor Luc Illusie}} 
\setlength{\baselineskip}{1.0\baselineskip}
\begin{abstract}
\noindent 
  We construct the fine moduli space of log abelian varieties, which gives a compactification of the moduli space of abelian varieties.
\end{abstract}
\section*{Contents}

\noindent \S \ref{sec:functor}. Moduli functors

\noindent \S \ref{s:pp}. Principal polarization

\noindent \S \ref{s:E(A)}. Universal additive extension

\noindent \S \ref{s:deformation}. Deformation of log abelian varieties

\noindent \S \ref{s:related}. Some related functors

\noindent \S \ref{s:Artin}. Artin's criterion

\noindent \S \ref{sec:pf-main}. Representability of moduli functors

\noindent \S \ref{s:val}. Valuative criterion 

\noindent \S \ref{s:toroidal}. Toroidal compactifications

\medskip

\section*{Introduction}
\renewcommand{\thefootnote}{\fnsymbol{footnote}}
\footnote[0]{Primary 14K10; 
Secondary 14J10, 14D06} 

  This part is the most important in our series of papers on log abelian varieties. 
  We construct the fine moduli space of principally polarized log abelian varieties with level structure and 
a prescribed admissible degeneration, 
and show that this moduli space is isomorphic
to a toroidal compactification of the moduli space of principally polarized abelian varieties with
level structure. 
  Our moduli space is obtained as a part of the moduli 
space without prescription of degeneration, which is a log algebraic space in the second sense.

  More precisely, we prove the following (1), (2), and (3).

(1) The moduli functor of $g$-dimensional principally polarized log abelian varieties with level $n$ structure is represented by a proper and log smooth log algebraic space over $\bZ[1/n]$ $(n \geq 3)$ in the second sense (Theorem \ref{thm:mainthm1}).

(2) The moduli functor of $g$-dimensional principally polarized log abelian varieties with level $n$ structure and with local monodromies in a prescribed admissible cone decomposition $\Sigma$ is represented by a proper and log smooth algebraic space with fs log structure over $\bZ[1/n]$ $(n \ge 3)$
(Theorem \ref{thm:mainthm2}).

(3) The space in (2) coincides with the toroidal compactification (\cite{FC}) endowed with a natural log structure of the moduli space of $g$-dimensional principally polarized abelian varieties with level $n$ structure associated with local monodromies in $\Sigma$ (Theorem \ref{thm:toroidal}).

  An important new instrument in this part of series of papers is the theory of deformations of log abelian varieties. 
  Like the other theories on log abelian varieties, this is also parallel to the corresponding classical theory of abelian varieties (see, for example, Proposition \ref{p:deform}), 
which reflects the remarkable feature of log abelian varieties that they are some kind of degenerations of abelian varieties but still behave as proper smooth group objects without degeneration.

  We are strongly influenced by the pioneering work by K.\ Fujiwara \cite{Fujiwara}.
  In fact, the prototypes of all our results and techniques of proofs can be found in his paper. 

  For other compactifications of the moduli space of abelian varieties, see \cite{Al}, \cite{Nakamura}, \cite{Olsson2}. 

  In Section \ref{sec:functor}, we describe main results. 
  In Section \ref{s:pp}, we show some results on principal polarizations. 
  In Section \ref{s:E(A)}, we introduce universal additive extensions of log abelian varieties and prove some basic facts.
  In Section \ref{s:deformation}, using them, we discuss the deformation of log abelian varieties. 
  After some preliminaries in Section \ref{s:related}, the core of the proof is an application of Artin's criterion to our moduli functor in Section \ref{s:Artin}. 
  We prove (1) and (2) in Section \ref{sec:pf-main} except properness of our moduli spaces. 
  In Section \ref{s:val}, we prove the properness.
  Finally, we prove (3) in Section \ref{s:toroidal}.

\smallskip

{\sc Acknowledgments.}
The first author is partially supported by JSPS, Kakenhi (C) No.\ 24540035, Kakenhi (C) No.\ 15K04811, and Kakenhi (C) No.\ 20K03555. 
The second author is partially  supported by NFS grants DMS 1303421 and DMS 1601861.  
The third author is partially supported by JSPS, Kakenhi (C) No.\ 22540011, 
Kakenhi (B) No.\ 23340008, and Kakenhi (C) No.\ 16K05093.

\section{Moduli functors}
\label{sec:functor}
  In this section, we state our main results. 

Fix an integer $g\ge 1$ and a finitely generated free abelian group $W$
of rank $g$. Let $S_{\Q}(W)$ be the $\Q$-vector space of symmetric
bilinear forms $W\times W\to \Q$.

\begin{defn}\label{defn:adm-decomp}
 (Cf.\ \cite[Chapter IV, Definition 2.2]{FC}.)
  An \emph{admissible cone decomposition} (or an \emph{admissible fan}) 
$\Sigma$ of
    $S_{\Q}(W)$ is a 
  set of finitely generated (sharp) $\Q$-cones in $S_{\Q}(W)$ satisfying the
  following conditions:

  (1) For $\sigma \in \Sigma$, every face of $\sigma$ is in $\Sigma$;

  (2) For $\sigma,\tau \in \Sigma$, the intersection $\sigma \cap
  \tau$ is a face of $\sigma$;

  (3) $\Sigma$ is stable under the action of $\Aut_{\Z}(W)$. Here
  $\alpha \in \Aut_{\Z}(W)$ acts on $S_{\Q}(W)$ by $b \mapsto
  b(\alpha(\cdot), \alpha(\cdot))$;

  (4) The number of the $\Aut_{\Z}(W)$-orbits in $\Sigma$ is finite;

  (5) For any $\sigma\in \Sigma$, any element of $\sigma$ is positive
  semi-definite, i.e., $b(w,w)\ge 0$ for any $b\in \sigma$ and any
  $w\in W$;

  (6) For each positive semi-definite symmetric bilinear form $b\colon
  W\times W \to \R$, there exists a unique $\sigma \in \Sigma$ for
  which $b$ is contained in the interior of $\sigma\otimes_{\Q_{\ge 0}}\R_{\ge 0}$. 
\end{defn}

\begin{sbrem}
(1) The last two conditions in Definition \ref{defn:adm-decomp}
can be replaced with the condition that 
the support $\bigcup_{\sigma \in \Sigma}\sigma$ of $\Sig$ is equal to the set of all positive semi-definite symmetric bilinear forms.

(2) It is known that an admissible cone decomposition of $W$ exists by
  the reduction theory (cf.\ \cite{Igusa}, \cite{Namikawa} 8.5).
\end{sbrem}

\begin{para}
  Recall that a polarization on a weak log abelian variety $A$ over an fs log scheme $S$ is a symmetric biextension of the pair $(A,A)$ by $\Gmlog$ 
whose pullback to $\ol s$ for any $s \in S$ is induced by a polarization on the log $1$-motif corresponding to $A\times_S\ol s$
(\cite{KKN5} 1.3).

  We say that a polarization is \emph{principal} if for any $s \in S$, 
it induces an isomorphism from $A \times_S\ol s$ to its dual.
  (The dual of a weak log abelian variety with constant degeneration is defined in \cite{KKN5} 1.2.)
\end{para}
 
We now define a principally polarized log abelian variety of degeneration along an admissible fan as follows.

\begin{defn}
  Let $\Sigma$ be an admissible cone decomposition of $S_{\Q}(W)$, and
  $A$ a principally polarized log abelian variety of dimension $g$
  over an fs log scheme $S$. We say that the \emph{local monodromies of $A$
  are in $\Sigma$}, or $A$ is \emph{compatible with $\Sigma$}, or $A$ is 
  \emph{of degeneration along $\Sigma$}, if strict \'etale locally on $S$, there exists a surjective
  homomorphism $f\colon W\to \ol Y$ satisfying the following
  condition: for any $s\in S$, there exists a $\sigma \in \Sigma$ such
  that, for any homomorphism $h\colon (M_S/\O_S^{\times})_{\ol
    s} \to \N$, the composition
  $$W\times W \overset{f\times f}{\to} \ol Y_{\ol s}\times 
\ol Y_{\ol
    s}\overset{\phi\times \id}{\to } \ol X_{\ol s}\times 
\ol Y_{\ol s}
  \overset{\langle\cdot,\cdot\rangle}{\to} (M_S^{\gp}/\O_S^{\times})_{\ol
    s}\overset{h^{\gp}}{\to} \Z$$
  is contained in $\sigma$. 
  Here $\langle\cdot,\cdot\rangle\colon \ol X \times \ol Y \to \Gmlog/\Gm$ is the canonical pairing of $\bZ$-modules determined by $A$ (\cite{KKN2} 4.4) and 
$\phi \colon \ol Y_{\ol s} \to
  \ol X_{\ol s}$ is the homomorphism defined by the polarization. 
\end{defn}

\begin{defn}
  Let $n$ be an integer $\ge 1$. 
  Let $S$ be an fs log scheme over $\Spec \bZ[1/n]$. 
  Let $A$ be a log abelian variety over $S$ of dimension $g$. 
  A \emph{level $n$ structure} of $A$ is an isomorphism $(\bZ/n\bZ)^{2g} \overset \cong \to \Ker(n\colon A \to A)$. 
\end{defn}

\begin{defn}
\label{defn:functor}
  Let $n$ be an integer $\ge 1$. We define a moduli functor
  $$F_{g,n}\colon (\fs/\Z[1/n])\to (\mathrm{set})$$ and a moduli
  functor $$F_{g,n,\Sigma}\colon (\fs/\Z[1/n])\to (\mathrm{set})$$ for
  an admissible cone decomposition $\Sigma$ of $S_{\Q}(W)$ as follows.
  For an object $U$ of $(\fs/\Z[1/n])$, we set
  \begin{align*}
    F_{g,n}(U) :=& \{g\text{-dimensional principally polarized log
      abelian variety over } U \\ & \text{ with level } n \text
    { structure }\}/\isom; \\
    F_{g,n,\Sigma}(U) := &\{g\text{-dimensional principally polarized
      log 
      abelian variety over } U \\ & \text{ with level } n \text
    { structure and with local monodromies in }\Sigma \}/\isom.
  \end{align*}
\end{defn}

  The following three theorems are the main results of this series of papers.

\begin{thm}\label{thm:mainthm1}
  If $n \ge3$, the moduli functor $F_{g,n}$ is represented by a proper and log smooth log algebraic space over $\bZ[1/n]$ in the second sense ({\rm{\cite{KKN4} 10.1}}). 
\end{thm}

  This was partially mentioned in \cite{KKN4} 10.6.
  Here a log algebraic space in the second sense is said to be {\it log smooth} if it admits a log smooth cover, i.e., we can take a log smooth $F'$ in the definition in \cite{KKN4} 10.1. 
  For the definition of the properness, see \cite{KKN4} 17.3.

\begin{thm}\label{thm:mainthm2}
  If $n \ge 3$, the moduli functor $F_{g,n,\Sigma}$ is represented by a proper and log smooth log algebraic space over $\Z[1/n]$ in the first sense ({\rm{\cite{KKN4} 10.1}}), that is, 
it is an algebraic space with fs log structure over $\Z[1/n]$, and it is proper and log smooth.
\end{thm}

\begin{thm}\label{thm:toroidal}
  If $n \ge 3$, the space representing the functor $F_{g,n,\Sigma}$ in Theorem $\ref{thm:mainthm2}$ coincides with the toroidal
  compactification $\ol {\mathcal A}_{g,n}$ associated to $\Sigma$ ({\rm{\cite{FC}}}) 
of the moduli space $\mathcal A_{g,n}$ of principally polarized abelian varieties with level $n$ structure, 
endowed with the fs log structure defined by the divisor $\ol {\mathcal A}_{g,n}-\mathcal A_{g,n}$. 
\end{thm}

\begin{para}
  We explain the relationship between these theorems. 
  Roughly speaking, Theorem \ref{thm:mainthm1} except properness is deduced from Theorem \ref{thm:mainthm2} except properness. 
  But the properness of $F_{g,n,\Sig}$ in Theorem \ref{thm:mainthm2} is deduced from that of $F_{g,n}$ in Theorem \ref{thm:mainthm1}.
  As for Theorems \ref{thm:mainthm2} and \ref{thm:toroidal}, 
it is possible to prove Theorem \ref{thm:toroidal} first using the theory of Faltings--Chai (\cite{FC}) and 
deduce Theorem \ref{thm:mainthm2} by checking some properties of $F_{g,n,\Sigma}$ (cf.\ the second last paragraph of \ref{1stpf-thm:toroidal}), as we did in the one-dimensional case \cite{KKN3}.
  But, below, we prove Theorem \ref{thm:mainthm2} first without the use of \cite{FC}, and deduce Theorem \ref{thm:toroidal} from it.
\end{para}

\begin{para}
  We hope that the theory of log abelian varieties works for variants such as PEL-type moduli problems, compactifications over $\bZ$, $X_0(N)$-type moduli problems, moduli stacks for $n=1, 2$, and so on. 
  In particular, in a forthcoming paper, we plan to throw the coefficient rings to the above three theorems so that the resulting space 
in the generalized version of Theorem \ref{thm:mainthm2} 
gives a moduli interpretation of the space constructed in 
\cite{Lan}.
  See also \cite{KKN3} Section 6.
\end{para}

\section{Principal polarization}
\label{s:pp}

  We give some propositions on principal polarizations, which will be used later. 

\begin{para}
\label{vertical}
  Let $S$ be an fs log scheme. 
  Consider $\Gmlog$ on $(\fs/S)_{\et}$. 
  Let $\Gmlog^{\mathrm{vert}}$ be the vertical part of it, that is, 
the subgroup sheaf of $\Gmlog$ consisting of the sections $x$ satisfying the following condition: There are $a, b \in M_S^{\gp}$ such that 
$a \vert x \vert b$, that is, both $a^{-1}x$ and $x^{-1}b$ belong to $M \subset \Gmlog$. 
\end{para}

\begin{lem}
\label{l:Gmlogvert}
  Let $A$ be a weak log abelian variety over an fs log scheme $S$. 
 
$(1)$ $\cH^0(A,\Gmlog/\Gmlog^{\mathrm{vert}})=\Gmlog/\Gmlog^{\mathrm{vert}}$.

$(2)$ $\cHom(A,\Gmlog/\Gmlog^{\mathrm{vert}})=0$. 
\end{lem}

\begin{pf}
  By \cite{KKN5} Lemma 3.1, (2) is reduced to (1). 

  We prove (1). 
  This is an analogue of a part of \cite{KKN5} Proposition 2.1 and proved similarly as follows. 

  First, to prove the case where $A$ is with constant degeneration, we use an analogue of \cite{KKN2} Proposition 7.9 (2). 
  That is, in the situation there, we have 
$$f_U^{-1}(M^{\gp}_{V(\Delta')_U}/M^{\gp,\mathrm{vert}}_{V(\Delta')_U}) \cong (M_U^{\gp}/\cO^{\times}_U \oplus \overline X)/({\rm the\ vertical\ part}),$$
where $U$ is any fs log scheme over $S$, $M^{\gp,\mathrm{vert}}_{V(\Delta')_U}$ is the restriction of $\Gmlog^{\mathrm{vert}}$ to the small \'etale site of $V(\Delta')_U$, \lq\lq(the vertical part)'' is the subsheaf of $M_U^{\gp}/\cO^{\times}_U \oplus \overline X$ consisting of the sections sent into $M^{\gp,\mathrm{vert}}_{V(\Delta')_U}$.
  But, since any interior element of $\cS$ is an interior element of 
$\Delta'{}\spcheck$ (still in the notation there), 
$\Delta'{}\spcheck \cdot \cS^{\gp}=C^{\gp}$ contains $X$ so that $X$ is always included in the vertical part. 
  Thus the right-hand-side coincides with $M_U^{\gp}/M_U^{\gp,\mathrm{vert}}$, where $M_U^{\gp,\mathrm{vert}}$ is the restriction of 
$\Gmlog^{\mathrm{vert}}$ to $U_{\et}$.
  Using this, in the notation of \cite{KKN5} 6.4, we can show 
$$g_*(\Gmlog/\Gmlog^{\mathrm{vert}}) = \Gmlog/\Gmlog^{\mathrm{vert}}.$$  
  Hence, for any $U$ over $S$, we have 
$$H^0(A_U,\Gmlog/\Gmlog^{\mathrm{vert}})=H^0(Y, H^0(B_U, \Gmlog/\Gmlog^{\mathrm{vert}}))=H^0(U,\Gmlog/\Gmlog^{\mathrm{vert}})$$
in the notation in \cite{KKN5} 6.4, which completes the proof of the case of constant degeneration. 

  In the general case, the zero section induces a surjection $H^0(A_U,\Gmlog/\Gmlog^{\mathrm{vert}})\to H^0(U,\Gmlog/\Gmlog^{\mathrm{vert}}).$
  To reduce the injectivity of this surjection to the case of constant degeneration, it is sufficient to show that $H^0(A_U,\Gmlog/\Gmlog^{\mathrm{vert}})\to \displaystyle\prod_{u \in U}H^0(A_u,\Gmlog/\Gmlog^{\mathrm{vert}})$ is injective. 
  We can replace $A_U$ here with any fs log scheme over it.  
  Then the injectivity follows. 
\end{pf}

\begin{rem}
  In the proof of the statement $\cH^0(A,\Gmlog/\Gm)=\Gmlog/\Gm$ which is a part of \cite{KKN5} Proposition 2.1, the corresponding 
reduction step to the case of the constant degeneration 
contains a rather inadequate explanation (cf.\ the first paragraph in \cite{KKN5} Proposition 12.1).  
  A correct argument is the one in the proof of the above lemma. 
\end{rem}

\begin{prop}
\label{p:ppfppf}
  Let $A$ be a weak log abelian variety over an fs log scheme $S$. 
  Any strict fppf locally given principal polarization with descent data gives a unique principal polarization on $A$. 
\end{prop}

\begin{pf}
  Note that this is easy if $\BIEXT(A, A;\Gmlog)$ is an fppf sheaf, but we have not yet proved that it is so.

  We reduce to the descent of an extension of $A$ by $\Gmlog$ as follows. 

  The data give a homomorphism fppf locally 
$$A \to \cExt(A,\Gmlog)$$
by \cite{KKN5} Proposition 2.3. 
  It is enough to show that the image $E$ of any section of $A$ by this homomorphism with the descent data can be uniquely descended. 

  First, we claim that $E$ as an extension comes from a unique extension $E'$ of $A$ by 
$\Gmlog^{\mathrm{vert}}$. 

  This claim is proved as follows. 
  Since 
$\cHom(A,\Gmlog/\Gm^{\mathrm{vert}})=0$ (Lemma \ref{l:Gmlogvert} (2)), we have the commutative diagram 
$$\begin{CD}
0 @>>>\cExt(A,\Gmlog^{\mathrm{vert}}) @>>>\cExt(A,\Gmlog)@>>> \cExt(A, \Gmlog/\Gmlog^{\mathrm{vert}})\\
@. @VVV @VVV @| \\
0 @>>> \cExt(A,\Gmlog^{\mathrm{vert}}/\Gm) @>>> \cExt(A,\Gmlog/\Gm) @>>>\cExt(A, \Gmlog/\Gmlog^{\mathrm{vert}})  \\
\end{CD}$$
with exact rows. %
  Hence, it is enough to show that 
the image of $E$ in $\cExt(A,\Gmlog/\Gm)$ comes from 
$\cExt(A,\Gmlog^{\mathrm{vert}}/\Gm)$.
  By \cite{KKN5} Proposition 1.6 (1) and \cite{KKN5} 2.7, this image comes from $\cHom(\overline Y, \Gmlog/\Gm)/\overline X$ by the composite 
$$\cHom(\overline Y, \Gmlog/\Gm)/\overline X \to \cExt(A/G,\Gmlog/\Gm) \to 
\cExt(A,\Gmlog/\Gm),$$
where the first arrow is the injection 
in \cite{KKN5} Lemma 2.6 (2). 
  Here the notation is as in there. 
  Further it comes from $\cHom(\overline Y, \Gmlog/\Gm)^{(\overline X)}/\overline X$ because it is so at each fiber. 

  Since the above composite is compatible with another composite 
$$\cHom(\overline Y, \Gmlog^{\mathrm{vert}}/\Gm)/\overline X \to \cExt(A/G,\Gmlog^{\mathrm{vert}}/\Gm) \to \cExt(A,\Gmlog^{\mathrm{vert}}/\Gm),$$
to complete the proof of the claim, it is sufficient to show 
$$\cHom(\overline Y, \Gmlog/\Gm)^{(\overline X)} \subset \cHom(\overline Y, \Gmlog^{\mathrm{vert}}/\Gm).$$
  Let $\varphi\colon \overline Y \to \Gmlog/\Gm$ be in the $(X)$-part. 
  Then for $y \in \overline Y$ in each fiber, there exist $x_1, x_2 \in \overline X$ such that $\langle x_1, y\rangle | \varphi(y) | \langle x_2, y\rangle$.
  Since $\langle x_i, y\rangle$ ($i=1,2$) are in the log of the base, this implies that $\varphi(y) \in \Gmlog^{\mathrm{vert}}$.
  Hence we see that $\varphi$ factors through $\Gmlog^{\mathrm{vert}}/\Gm$, as desired.

  Let $\tilde E'$ be the extension of $\tilde A$ by $\Gmlog^{\mathrm{vert}}$ induced by our $E'$.  
  Since the composite 
$$\cHom(\overline Y, \Gmlog^{\mathrm{vert}}/\Gm)/\overline X \to \cExt(A/G,\Gmlog^{\mathrm{vert}}/\Gm) \to \cExt(\widetilde A/G,\Gmlog^{\mathrm{vert}}/\Gm)$$
is zero, $\tilde E'$ induces the trivial extension of $\tilde A$ by $\Gmlog^{\mathrm{vert}}/\Gm$. 
  Fix a splitting of this trivial extension and identify the total space of this extension with the product $\Gmlog^{\mathrm{vert}}/\Gm \times \tilde A$. 

  It is enough to descend $E'$. 
  We can do it similarly as in the case of $A$ itself (cf.\ \cite{KKN4} Section 9).  
  We may assume that a prime $\ell$ is invertible on the base. 
  We cover $\tilde E'$ with respect to the k\'et topology by the images with respect to the multiplication by $\ell^n$ ($n\ge0$) of a representable object explained below. 
  Then as in the same way as in the case of $A$ (cf.\ \cite{KKN4} Section 9), we descend $\tilde E'$ by descending this representable object with a partial group law, we descend the kernel of $\tilde E' \to E'$, and we divide $\tilde E'$ by the kernel to obtain the descended $E'$.
  The representable object we use is as follows. 
  Let $I$ be a model of $\tilde A$ with respect to a wide cone (as in \cite{KKN4} Section 9).
  Let $a$ be an interior in the log of the base.
  Let $L$ be the part of $\Gmlog^{\mathrm{vert}}/\Gm$ consisting of $x$ such that $a^{-1} \vert x \vert a$. 
  Consider the inverse image in $\tilde E'$ of $L \times I$.
  Then it is representable. 
  We explain the representability. 
  Assume that $A=S$. 
  Then the inverse image of $L$ in $\Gmlog^{\mathrm{vert}}$ is represented by the fs log scheme defined as the fiber product of 
$$S \to \Spec(\bZ[\bN]) \leftarrow \Spec(\bZ[\bN^2]),$$ 
where the first morphism sends $1 \in \bN$ to $a^2$ in $M_S$, and the second morphism is induced by the diagonal homomorphism $\bN \to \bN^2$. 
  The representability in the general case follows from this.
\end{pf}

\begin{rem}
  (1) The final portion of the above proof can be regarded as a generalization of the argument with models for a log abelian variety in \cite{KKN4} to that for a kind of log semiabelian variety. 
  Though we even have not yet defined a log semiabelian variety, the theory of (proper-like or quasiprojective) models of a log semiabelian variety should be an important subject to study. 

  (2) In \cite{KKN4} Section 9, Lemma 9.10 and Lemma 9.11 are not essential. 
  In fact, in \cite{KKN4} 9.12, we may assume that the conditions (5)--(7) are satisfied only by replacing the index $\lambda$. 
\end{rem}

\begin{prop}
\label{p:pplimit}
  Let $(S_{\lambda})_{\lambda}$ be a filtered projective system of quasicompact and quasiseparated 
fs log schemes whose transition morphisms 
are affine and strict.
  Let $S:=\varprojlim S_{\lambda}$.
  For a principally polarized log abelian variety $(A,p)$ over $S$, there are an index $\lambda$ and a principally polarized log abelian variety $(A_{\lambda},p_{\lambda})$ over $S_{\lambda}$ whose pullback to $S$ is isomorphic to $(A,p)$. 

  Further, for some $\lambda$ and two 
principally polarized log abelian varieties over $S_{\lambda}$ whose pullbacks to $S$ are isomorphic to each other, there is another index $\lambda \to \lambda'$ such that their pullbacks to $S_{\lambda'}$ are already isomorphic to each other. 
\end{prop}

\begin{pf}
  By \cite{KKN4} Proposition 9.2, we may assume that $A$ over $S$ comes from an essentially unique weak log abelian variety $A_{\lambda}$ over some $S_{\lambda}$.
  The problem is to spread out $p$. 
  (We mean that there exist an essentially unique principal polarization on $A_{\lambda}$ for a sufficiently large $\lambda$ whose pullback coincides with $p$.)
  This is shown by the same method in the proof of the previous proposition. 

  First, we can replace $\Gmlog$ by $\Gmlog^{\mathrm{vert}}$ in $\Biext(A,A;\Gmlog)$.
  In fact, by the claim in the proof of Proposition \ref{p:ppfppf}, the polarization lives in $\Hom(A,\cExt(A,\Gmlog^{\mathrm{vert}}))$.
  Since $\Hom(A,\Gmlog^{\mathrm{vert}}) \subset \Hom(A,\Gmlog)=0$ (\cite{KKN5} Proposition 2.1 and \cite{KKN5} Lemma 3.1), we have 
$\Hom(A,\cExt(A,\Gmlog^{\mathrm{vert}}))=\Biext(A,A;\Gmlog^{\mathrm{vert}})$ by \cite{KKN5} Lemma 3.4.

  Second, we consider the representable object induced by $L \times I \times I$ in the notation in the proof of Proposition \ref{p:ppfppf}. 
  Then we can spread out the biextension with partial group laws in the same way. 
  To prove that the biextension on $S_{\lambda}$ is a principal polarization for a sufficiently large $\lambda$, 
we may assume that the biextension on $S_{\lambda}$ is symmetric. 
  Further, we may assume that there is a chart of $S_{\lambda}$. 
  Then, by dividing $S_{\lambda}$ into a finite number of constant log loci, we may assume that $S_{\lambda}$ has the constant log. 
  Hence we can replace the log abelian variety concerned with the corresponding log $1$-motif $[Y \to G_{\log}]$.
  Then, it is enough to show that the corresponding homomorphism $[Y \to G_{\log}] \to [Y \to G_{\log}]$ is an isomorphism over $S_{\lambda}$ for a sufficiently large $\lambda$ when its pullback to $S$ is an isomorphism. 
  It is valid because $Y$ is a locally constant sheaf and 
the homomorphism $G_{\log} \to G_{\log}$ comes from a homomorphism $G \to G$ (\cite{KKN2} Proposition 2.5). 
\end{pf}

\begin{rem}
  We have not yet had the property of local finite presentation for nonk\'et $\cB iext$.
  If it is the case, we can use it. 
  Or, if we had the general theory of duals for a weak log abelian variety not necessarily with constant degeneration, the proposition follows from \cite{KKN4} Proposition 9.2 (1). 
\end{rem}

\begin{prop}
\label{p:ppinj}
  Let $A$ be a weak log abelian variety over an fs log scheme $S$. 
  Then the homomorphism $A \to \cExt(A,\Gmlog)$ induced by a principal polarization on $A$ ({\rm{\cite{KKN5}}} Proposition $2.3$) is injective.
\end{prop}

  To prove this, we need a lemma, which also will be used in later sections. 

\begin{lem}
\label{l:endo}
  Let $S$ be an fs log scheme, $A_1$ and $A_2$ weak log abelian varieties over $S$ with constant degeneration, $f$ and $g$ homomorphisms from $A_1$ to $A_2$.  
  If $f_{\overline s}=g_{\overline s}$ for any $s \in S$, 
then $f=g$. 
\end{lem}

\begin{pf}
  Let $M_i=[Y_i \to G_{i,\log}]$ be a log 1-motif over $S$ corresponding to $A_i$ for $i=1,2$. 
  It suffices to prove that a homomorphism $M_1 \to M_2$ whose pullback to $\overline s$ is zero for any $s \in S$ is zero.
  The homomorphism consists of a homomorphism $Y_1 \to Y_2$ and that of $G_{1,\log}\to G_{2,\log}$. 
  The former is zero because $Y_i$ are locally constant for $i=1,2$. 
  The latter is induced by a homomorphism $G_1 \to G_2$ by \cite{KKN2} Proposition 2.5, 
which induces a homomorphism of the torus parts and the abelian parts.
  The homomorphism of the torus parts is zero because the corresponding homomorphism of the character groups is so.
  The homomorphism of the abelian parts is also zero. 
  Since there is no nontrivial homomorphism from the abelian part of $G_1$ to the torus part of $G_2$, we conclude that the homomorphism $G_1 \to G_2$ is zero. 
  Hence the homomorphism of $M_1 \to M_2$ is zero.
\end{pf}

\begin{para}
  We prove Proposition \ref{p:ppinj}. 
  By Proposition \ref{p:pplimit}, we may assume that the base is finitely generated over $\bZ$. 

  It is enough to prove that for any fs log scheme $U$ over the base, the homomorphism $A(U) \to \Ext(A_U,\Gmlog)$ is injective.

  First, $\cExt(A,\Gmlog)$ is locally of finite type. 
  This is by \cite{KKN5} Proposition 12.8 (1) and the fact that $\cExt(A,\Gmlog) \to \cExt_{\ket}(A,\Gmlog)$ is injective. 
  The last fact is seen as follows. 
  Since both $A$ and $\Gmlog$ are k\'et sheaves (\cite{KKN4} Theorem 5.1 and \cite{Kato:FI2} Theorem 3.2), any extension on the \'etale site is a k\'et sheaf, too. 
  Hence the splitting on the k\'et site gives a splitting on the \'etale site. 

  On the other hand, $A$ is locally of finite presentation by \cite{KKN5} Proposition 12.7.
  Hence, we may assume that $U$ is finitely generated over $\bZ$. 

  Since $A$ is locally of finite type, we can replace $U$ with the strict localization at each point.
  Again by that $A$ is locally of finite type and by Artin's approximation theorem, we can replace $U$ with a complete noetherian local ring. 
  Here we have the GAGF for $H^0(U,A)$, that is the algebraization of the sections, which can be proved by taking a model.
  By it, we further replace $U$ with an Artin local ring so that we may assume that $A$ is with constant degeneration (\cite{KKN2} Theorem 4.6 (2)). 

  Thus, it is enough to show the proposition under the additional assumption that $A$ is with constant degeneration. 
  We denote by $A^*$ the dual of $A$. 
  In this case, $A \to \cExt(A,\Gmlog)$ factors through $f\colon A \to A^*$, where $A^*$ is identified with a subgroup of $\cExt(A,\Gmlog)$ by \cite{KKN2} Remark 7.5 (2). 
  By the assumption, $f_{\overline s}$ is an isomorphism for any $s \in S$. 
  Then, $f$ is an isomorphism by Lemma \ref{l:endo} as desired. 
\end{para}

\section{Universal additive extension}
\label{s:E(A)}
  We introduce the universal additive extension of a log abelian variety and will use it for the study of the deformation of log abelian variety in the next section.
  See \cite{Mochizuki} for a log elliptic curve case (though the definition of log elliptic curve in \cite{Mochizuki} is different from ours).

\begin{para}
  We define the Lie sheaf $\Lie(A)$ and coLie sheaf $\coLie(A)$.
  Let $S$ be an fs log scheme. 

  In general, for a group sheaf $G$ on $(\fs/S)_{\et}$, we define the Lie sheaf $\Lie(G)$ on $(\fs/S)_{\et}$ by 
$$\Lie(G)(U) = \Ker(G(U[\varepsilon]/(\varepsilon^2)) \to G(U)),$$
where $U[\varepsilon]/(\varepsilon^2)$ denotes the fs log scheme whose underlying space is that of $U$, whose structure sheaf is $\cO_U[\varepsilon]/(\varepsilon^2)$ and whose log structure is the pullback of that of $U$. 

  We have a natural action of $\O_S$ on $\Lie(G)$ (cf.\ \cite{KKN1} 1.3.11).  
  In many cases, $\Lie(G)$ is an $\cO_S$-module via this action. 
  For example, we have $\Lie(\Gmlog)=\Lie(\Gm)=\Ga$. 
\end{para}

\begin{para}
  If $A$ is a weak log abelian variety over $S$, $\Lie(A)$ is a locally free $\cO_S$-module of rank $\dim(A)$ because it coincides with the Lie sheaf of the semiabelian part.

  Let $$\coLie(A)=\cHom_{\cO_S}(\Lie(A),\cO_S),$$ which is also a locally free $\cO_S$-module of rank $\dim(A)$.
\end{para}

\begin{para}
  We introduce the dual of $A$ in case where $A$ is principally polarized. 
  In the rest of this section, we assume that $A$ is principally polarized unless explicitly stated otherwise. 
  We denote by $A^*$ the image of $A \to \cExt(A,\Gmlog)$ and call it the {\it dual} of $A$. 
  By Proposition \ref{p:ppinj}, $A^*$ is isomorphic to $A$.  
  We identify $A^*$ with $A$ via this isomorphism.
\end{para}

\begin{prop}
\label{p:coLie}
  Let $A$ be a principally polarized log abelian variety over an fs log scheme $S$. 
  Then 
$$\coLie(A^*) \cong \cExt(A,\Ga)^*:=\cHom(\cExt(A,\Ga),\Ga)$$ 
and it is a vector bundle of rank $\dim A$. 
\end{prop}

\begin{pf}
  The map $A^* \to \cExt(A,\Gmlog)$ induces a map 
$$\Lie(A^*) \to \Lie(\cExt(A,\Gmlog)) \to \cExt(A,\Lie(\Gmlog)) = \cExt(A,\Ga).$$ 
  It is enough to show that it is an isomorphism. 

  First we prove it under the assumption that $A$ is with constant degeneration. 
  By \cite{KKN2} Theorem 7.3 (1), the exact sequence 
(the notation is as usual as in there)
\begin{equation*}\tag{1}  
0 \to T_{\log}^{(Y)} \to \tilde A \to B \to 0 
\end{equation*}
yields isomorphisms $\cHom(\tilde A, \Ga)=0$ and $\cExt(\tilde A, \Ga)=\cExt(B, \Ga)$.
  Together with an exact sequence 
\begin{equation*}\tag{2}  
0 \to Y \to \tilde A \to A \to 0,
\end{equation*}
we have an exact sequence
\begin{equation*}\tag{3}  
0\to \cHom(Y, \Ga) \to \cExt(A,\Ga) \to \cExt(B, \Ga) \to 0.  
\end{equation*}
  Hence $\cExt(A,\Ga)$ is a vector bundle of the desired rank. 

  Further, from (2), we have $\Lie(A)= \Lie(\tilde A)$. 
  Together with (1), we have an exact sequence 
$$0 \to \Lie(T) \to \Lie(A) \to \Lie(B) \to 0.$$  
  Hence we have a commutative diagram 
$$\begin{CD}
0 @>>> \Lie(T^*) @>>> \Lie(A^*) @>>> \Lie(B^*) @>>> 0 \\
@. @VVV @VVV @VVV \\
0 @>>> \cHom(Y,\Ga) @>>> \cExt(A,\Ga) @>>> \cExt(B,\Ga) @>>> 0 
\end{CD}$$
with exact rows. 
  Since the left and the right vertical arrows are isomorphisms, the middle one is also, which completes the case of constant degeneration. 

  To prove the general case, by Proposition \ref{p:pplimit}, we may assume that $S$ is of finite type over $\bZ$.
  It suffices to show that
  $\Lie(A^*)(U) \to \Ext(A_U,\Ga)$ is bijective
for any fs log scheme $U$ over $S$.  

  We may assume that we are given the usual data by which we can discuss models. 
  (See \cite{KKN5} 4.5. 
  It is always the case after going strict \'etale locally over the base $S$.) 
  We prove the case where $U$ is the spectrum of a complete noetherian local ring $(R,m)$. 
  Consider the diagram 
$$\begin{CD}
  \Lie(A^*)(U) @>>> \Ext(A_U,\Ga) \\
@| @VVV \\
\varprojlim \Lie(A^*)(U_n) @>>> \varprojlim \Ext(A_{U_n},\Ga), 
\end{CD}$$
where $U_n$ is the spectrum of $R/m^{n+1}$ endowed with the pullback log structure from $U$.
  Since we already proved the case of constant degeneration, the bottom arrow is bijective.
  We prove that the right vertical arrow is injective so that bijective.
  To see it, since $\Ext(A_U,\Ga)$ is a subgroup of $H^1(A_U,\Ga)$ (in fact, the former coincides with the latter; see Lemma \ref{l:cubic}), it is enough to have GAGF for $H^1(A,\Ga)$.
  We take a complete and wide model $P$ of $A_U$ (\cite{KKN5} Proposition 10.3). 
  We take a prime $\ell$ which is invertible on $U$.
  For $i \ge0$, let $P_i \to A_U$ be the composite $P \to A_U \overset {\ell^i} \to A_U$. 
  Cover $A_U$ with respect to the k\'et topology by the disjoint union $\coprod P_i$ of the $P_i$.
  Consider the \v Cech-derived spectral sequence for $\coprod P_i/A_U$ 
$$E_2^{pq}=\check H^p({\textstyle\coprod} P_i/A_U, \cH^q(-,\Ga)) \Rightarrow H^{p+q}(A_U, \Ga).$$
  Note that the k\'et cohomology here coincides with the \'etale cohomology (cf.\ \cite{IKN} Proposition 3.7). 
  Since $E_1^{pq}$ is a projective limit of modules of finite lengths for every $p, q$, the classical GAGF implies that for $H^m(A_U,\Ga)$ for every $m$, which completes the proof of the current case.

  We prove the case of general $U$.
  We can prove that $\cExt(A,\Ga)=\cExt_{\ket}(A,\Ga)$ is locally of finite presentation by the method of \cite{KKN5} Proposition 12.8 (1). 
  Hence we may assume that $U$ is a strict localization of an fs log scheme of finite type over $\bZ$. 
  Then, 
again by that $\cExt(A,\Ga)$ is locally of finite presentation and 
by Artin's approximation theorem, we reduce to the previous case where $U$ is the spectrum of a complete noetherian local ring. 
\end{pf}

\begin{para}
  By Proposition \ref{p:coLie}, we have a canonical isomorphism 
$$\cExt(A,\cExt(A,\Ga)^*)=\cExt(A,\Ga) \otimes \cExt(A,\Ga)^*=
\cHom(\cExt(A,\Ga),\cExt(A,\Ga)).$$

  Define
$$0 \to \cExt(A,\Ga)^* \to E(A) \to A \to 0$$
as the extension associated with the canonical element in the left-hand-side 
corresponding to the identity in the right-hand-side of the above isomorphism.
  We call this the {\it universal additive extension} of $A$. 

  Then, by construction, it has a universality, that is, any extension of $A$ by some vector bundle $V$ of finite rank can be obtained by the pushout with respect to a unique homomorphism 
$\cExt(A,\Ga)^* \to V$. 
\end{para}

\begin{para}
  By Proposition \ref{p:coLie} and the identification $A=A^*$, the universal additive extension is isomorphic to 
$$0 \to \coLie(A) \to E(A) \to A \to 0.$$
\end{para}

  In the rest of this section, we prove the following proposition. 

\begin{prop}
\label{p:E(A)unique}
  Let $S \to S'$ be a strict nil immersion of affine fs log schemes of finite type over $\bZ$ by an ideal $I$ such that $I^2=0$. 
  Let $A$ be a principally polarized log abelian variety over $S$.
  Let $A'$ be a lift of $A$, that is, a principally polarized log abelian variety over $S'$ endowed with an isomorphism between $A$ and its pullback to $S$.  
  Then we have a natural isomorphism between $E(A')$ and the sheaf of groups $\cExt^{\natural}(A,\Gmlog)'$ of $\natural$-extensions defined below.
  In particular, the group sheaf $E(A')$ is independent of a lift up to isomorphisms. 
\end{prop}

\begin{para}
\label{Extnatural}
  We define $\cExt^{\natural}(A,\Gmlog)'$. 
  Consider the sheaf $\cH^1_{\mathrm{logcrys}}(A,\Gmlog)$ 
associated to the presheaf 
$(\fs/S') \ni U' \mapsto 
H^1(((\fs/A_{U'\times_{S'}S})/U')^{\log}_{\mathrm{crys}}, \Gmlog)$
(cf.\ \cite{Kato:FI1} Section 5 for the log crystalline site).  
  Let $\cExt^{\natural}(A,\Gmlog)$ be the fiber product of 
$$\cH^1_{\mathrm{logcrys}}(A,\Gmlog) \to 
\cH^1(A,\Gmlog) \leftarrow A.$$

  Consider the projection $q\colon 
(\fs/S')_{\et} \to 
(\fs/S')'_{\et}$, where the target is the subsite of $(\fs/S')_{\et}$ 
consisting of fs log schemes over $S'$ which are of finite type over $S'$. 
  Let $\cExt^{\natural}(A,\Gmlog)'=q^{-1}q_*\cExt^{\natural}(A,\Gmlog)$. 

  Further, in the rest of this section, for an fs log scheme $S$ of finite type over $\bZ$ and a sheaf $F$ on $(\fs/S)_{\et}$, we denote $q^{-1}q_*F$ by $F'$, where 
$q$ is the projection from $(\fs/S)_{\et}$ to its subsite consisting of fs log schemes over $S$ which are of finite type over $S$.
  For all $F'$ appearing in the rest of this section, it is plausible that $F = F'$. 
  In particular, it is plausible $\cExt^{\natural}(A,\Gmlog)' =\cExt^{\natural}(A,\Gmlog)$.
\end{para}

\begin{para}
  To prove Proposition \ref{p:E(A)unique}, we use another construction of $E(A)$ (see Proposition \ref{ue2}). 
  To explain this, we have to introduce the small site of $A$. 

  Let $A$ be a weak log abelian variety over an fs log scheme $S$. 

  We define the small site $(\llet/A)_{\et} \subset (\fs/A)_{\et}$ 
of $A$ as the full subcategory consisting of all log \'etale objects. 
  Here $(V \to A) \in (\fs/A)$ is said to be {\it log \'etale} if for any fs log scheme $W$ and any morphism $W \to A$, the fiber product $W \times_AV$ is a log \'etale log algebraic space over $W$ in the first sense (\cite{KKN4} 10.1).

  Next we define the site $(\fs\llet/A)_{\et}$ as follows. 
  Objects are the pairs $(U,V)$, where $U$ is an object of $(\fs/S)$ and $V$ is an object of $(\llet/A_U)$, where $A_U = A\times_SU$.
  Morphisms from $(U',V')$ to $(U,V)$ is a compatible pair of morphisms $U' \to U$ and $V' \to V$. 
  Coverings are the $((U_i,V_i)\to (U,V))_i$ such that every $U_i \to U$ is strict \'etale and that $(V_i \to V)_i$ is a strict \'etale covering. 
  
  Assume that we are given the usual data by which we can discuss models. 
  Then any log \'etale fs log scheme over some model of $A$ is an object of $(\llet/A)_{\et}$.
  This is because for a model $P$, the projection $P \to A$ is relatively represented by log blow-ups. 
  Further, under the condition in \cite{KKN5} 1.4.1, the proof of \cite{KKN5} Proposition 11.1 gives 
a set of topological generators consisting of these objects together with log \'etale fs log schemes over models of weak log abelian varieties which are isogeneous to $A$. 
  Similarly, based on \cite{KKN5} Proposition 11.1,  we can give a set of generators of $(\fs\llet/A)_{\et}^{\sim}$ consisting of the set of the pairs $(U,V)$, where $U$ is an object of $(\fs/S)$ and $V$ is a log \'etale fs log scheme over some model of a weak log abelian variety being isogeneous to $A_U$. 

  We can develop the theory of quasicoherent sheaves on $(\fs\llet/A)$.
  Let $\cO_A$ be the sheaf on $(\fs\llet/A)_{\et}$ defined by $\cO_A((U,V))=\Gamma(V,\cO_V)$. 
  Define a sheaf of $\cO_A$-modules $\omega^1_{A/S}$ by $\omega^1_{A/S}((U,V))=\Gamma(V,\omega^1_{V/U})$. 
  The sheaves $\Ga, \Gm, \Gmlog$ have been considered on $(\fs/A)_{\et}$ in former parts in our series of papers. 
  But they are also defined on $(\fs\llet/A)$ (by abuse of notation) by $\Ga=\cO_A, \Gm((U,V))=\Gamma(V,\cO_V^{\times})$, and 
$\Gmlog((U,V))=\Gamma(V,M^{\gp}_V)$.

  Let $\cH^q(A, -)$ be the $q$-th derived functor of the direct image functor 
$$(\fs\llet/A)_{\et} \to (\fs/S); F \mapsto (U \mapsto \varprojlim_{V/A_U}F((U,V))).$$
  Note that $\cH^q(A, \Ga)$, $\cH^q(A, \Gm)$, and $\cH^q(A, \Gmlog)$ in this sense coincide with the ones in the former parts of this series of papers, where the cohomology is considered via $(\fs/A)_{\et}$.
  This is because the canonical projection $(\fs/A)_{\et}^{\sim} \to (\llet/A)_{\et}^{\sim}$ preserves $\Ga$, $\Gm$, and $\Gmlog$. 
  Note that we do not consider $\omega^1_{A/S}$ on the big site $(\fs/A)_{\et}$. 

  Below, we also consider the small site of $\widetilde A$, the sheaf $\omega^1_{\widetilde A/S}$ on it, and its cohomologies, which are defined similarly. 
\end{para}

\begin{prop}
\label{p:lacoh2}
  Let $S$ be an fs log scheme which is of finite type over $\bZ$. 
  Let $A$ be a principally polarized log abelian variety over an fs log scheme $S$. 
  Then we have $\cH^0(A,\omega^1_{A/S})'\cong\coLie(A)$. 
  (See $\ref{Extnatural}$ for the notation $(-)'$.)
\end{prop}

\begin{para}
\label{pf-p:lacoh2}
  We prove this proposition till \ref{endpf_lacoh2}.
  Consider the natural map 
$$\cH^0(A,\omega^1_{A/S}) \to \cH^0(G,\omega^1_{G/S}) \to \coLie(G) \cong \coLie(A),$$
where $G$ is the semiabelian part of $A$. 

  First we assume that $A$ is of constant degeneration and we prove that the above map is bijective so that we also have $\cH^0(A,\omega^1_{A/S})'\cong\coLie(A)$.
  (Note that in this case, we do not use the assumption that $S$ is of finite type over $\bZ$.) 
  By \cite{KKN5} 6.4.1, we have an exact sequence 
$$0\to \cH^0(B, \omega^1_{B/S}) \to \cH^0(\tilde A, \omega^1_{\tilde A/S}) \to \cO_{S} \otimes_{\Z} X \to \cH^1(B,\omega^1_{B/S})$$
(the notation is as in there).
  Here the last arrow is zero because it factors through the abelian variety $\cExt(B,\Gm)$ and $\cH^1(B,\omega^1_{B/S})$ is a vector bundle.
  Thus we obtain the upper row of the commutative diagram 
$$\begin{CD}
0 @>>> \cH^0(B, \omega^1_{B/S}) @>>> \cH^0(\tilde A, \omega^1_{\tilde A/S}) @>>> \cO_S \otimes_{\Z} X @>>> 0 \\
@.@VVV @VVV @| @. \\
0 @>>> \coLie(B) @>>> \coLie(\tilde A) @>>> \cO_S \otimes_{\Z} X @>>> 0. \\
\end{CD}
$$
  Since the left vertical arrow is bijective, the middle one is also.

  Now we prove that $Y$ acts on $\cH^0(\tilde A, \omega^1_{\tilde A/S})$ trivially so that we may assume that the upper row splits.
  Once we have it, since $\cH^0(A, \omega^1_{A/S})=H^0(Y, \cH^0(\tilde A, \omega^1_{\tilde A/S}))$, we will have 
$\cH^0(A, \omega^1_{A/S})=\cH^0(\tilde A, \omega^1_{\tilde A/S})=\coLie(\tilde A)=\coLie(A)$, as desired. 
  Since $\cH^0(\tilde A \times \tilde A, \omega^1_{\tilde A\times \tilde A/S}) \to \coLie(\tilde A\times \tilde A)$ is injective as seen as in the same way, and since there is a natural decomposition $\coLie(\tilde A \times \tilde A)=\coLie(\tilde A) \times \coLie(\tilde A)$, we have the natural decomposition 
$$\cH^0(\tilde A \times \tilde A, \omega^1_{\tilde A\times \tilde A/S})=\cH^0(\tilde A,\omega^1_{\tilde A/S}) \times\cH^0(\tilde A,\omega^1_{\tilde A/S}).$$
  Let $y \in Y$ and $c_y\colon \tilde A \to \tilde A; a \mapsto y$ the constant map. 
  The translation by $y$ on $\cH^0(\tilde A,\omega^1_{\tilde A/S})$ coincides with $\cH^0(\cdot,\omega_{(\cdot)})$ of $\tilde A \overset{(\mathrm{id},c_y)}\to 
\tilde A \times \tilde A \overset{\mathrm{sum}}\to \tilde A$, which factors, via the above decomposition, through $\cH^0(\tilde A,\omega_{\tilde A/S}) \times \cH^0(\tilde A,\omega_{\tilde A/S}) \to \cH^0(\tilde A,\omega_{\tilde A/S}); \alpha \mapsto \alpha + (c_y^*\alpha)$. 
  Since $c_y$ factors through the base, $c_y^*\alpha$ is zero. 
  Thus $Y$ acts on $\cH^0(\tilde A, \omega^1_{\tilde A/S})$ trivially, which 
completes the proof of the case of constant degeneration. 
\end{para}

\begin{para}
\label{inj}
  To prove the general case of Proposition \ref{p:lacoh2}, 
  it is enough to show that 
$$(*) \qquad H^0(A,\omega^1_{A/S}) \to \coLie(A)(S)$$ 
is bijective.
  We may assume that we are given the usual data by which we can discuss models. 
  We prove the injectivity. %
  The GAGF for $H^0(A,\omega^1_{A/S})$ is proved in the same way as the GAGF for $H^1(A,\Ga)$ in the proof of Proposition \ref{p:coLie}. 
  By this and the constant degeneration case proved in \ref{pf-p:lacoh2}, we see that 

\smallskip

(1) $H^0(A_{\hat S},\omega^1_{A_{\hat S}/{\hat S}}) \to \coLie(A)(\hat S)$ is bijective, where $\hat S$ is the completion of the strict localization of $S$ at every point. 

\smallskip

  We prove that 

\smallskip

(2) $H^0(A_{S^{\mathrm{sh}}},\omega^1_{A_{S^{\mathrm{sh}}}/{S^{\mathrm{sh}}}}) \to \coLie(A)(S^{\mathrm{sh}})$ is injective, where $S^{\mathrm{sh}}$ is the strict localization of $S$ at every point. 

\smallskip

\noindent  Let $\alpha \in H^0(A_{S^{\mathrm{sh}}},\omega^1_{A_{S^{\mathrm{sh}}}/{S^{\mathrm{sh}}}})$. 
  Assume that the image of $\alpha$ in $\coLie(A)(S^{\mathrm{sh}})$ is zero. 
  Then, by (1), the image of $\alpha$ in $H^0(A_{\hat S},\omega^1_{A_{\hat S}/\hat S})$ is zero. 
  Take any fs log scheme $P$ which is of finite type over $S^{\mathrm{sh}}$ and any morphism $P \to A_{S^{\mathrm{sh}}}$ over $S^{\mathrm{sh}}$.
  The image of $\alpha$ in $H^0(P_{\hat S},\omega^1_{P_{\hat S}/{\hat S}})$ is clearly zero. 
  Then by Artin's approximation theorem, we see that 
$H^0(P,\omega^1_{P/{S^{\mathrm{sh}}}})$ is zero. 
  Since when $P$ varies, they cover $A_{S^{\mathrm{sh}}}$ (see \cite{KKN5} Proposition 11.1; actually, the proper $P$s already cover $A_{S^{\mathrm{sh}}}$), this implies that $\alpha$ is zero. 
  Thus (2) follows.

  Similarly, we can deduce from (2) the desired injectivity as follows. 
  Let $\alpha \in H^0(A, \omega^1_{A/S})$.
  Assume that the image of it in $\coLie(A)(S)$ is zero. 
  Then, by (2), the image of it in $H^0(A_{S^{\mathrm{sh}}},\omega^1_{A_{S^{\mathrm{sh}}}/{S^{\mathrm{sh}}}})$ is zero. 
  Take any fs log scheme $P$ which is of finite type over $S$ and any morphism $P \to A$ over $S$.
  The image of $\alpha$ in $H^0(P_{S^{\mathrm{sh}}},\omega^1_{P_{S^{\mathrm{sh}}}/{S^{\mathrm{sh}}}})$ is clearly zero. 
  Then the image of $\alpha$ in $H^0(P,\omega^1_{P/S})$ is \'etale locally vanishes, and hence, vanishes. 
  Since $P$s cover $A$, this implies that $\alpha$ is zero. 
\end{para}

  To prove the surjectivity of $(*)$, we use the following lemma. 

\begin{lem}
\label{l:flat}
  Let $A$ be a principally polarized log abelian variety over an fs log scheme $S$. 
  Then strict \'etale locally on $S$, there is a log regular fs log scheme $S_0$ of finite type over $\bZ$, a strict morphism $S \to S_0$, and 
a principally polarized log abelian variety $A_0$ over $S_0$ such that $A$ is isomorphic to $A_0\times_{S_0}S$.
\end{lem}

\begin{pf}
  By Proposition \ref{p:pplimit}, we may assume that $S$ is of finite type over $\bZ$. 
  We may work strict \'etale locally on $S$. 
  Taking a chart of $S$ and an admissible pairing valued in the groupification of the chart which induces the admissible pairing $\langle\cdot,\cdot\rangle$ associated to $A$, 
we may assume that there is a toric variety $S_0$ over $\bZ$ and a strict morphism $S \to S_0$ such that $\langle\cdot,\cdot\rangle$ comes from $S_0$. 
  We may further assume that $S \to S_0$ is a strict closed immersion. 

  For each point $s$ of $S$, consider the strict localizations at $s$  and their completions of both $S$ and $S_0$. 
  The $A$ lifts formally over the completion of the strict localization $S_0^{\mathrm{sh}}$ of $S_0$ by the formal smoothness of the local moduli in \cite{KKN6} Theorem 2.4. 
  By GAGF (\cite{KKN6} Theorem 6.1 and \cite{KKN6} Remark 6.1.1), it lifts also algebraically.  
  Then, by Proposition \ref{p:pplimit}, $A$ lifts over a subring of the completion which is finitely generated over $S_0^{\mathrm{sh}}$. 
  By Artin's approximation theorem, we take compatible sections over $S$ and over $S_0$, pull back the lift, and we have a lift of 
$A$ over $S_0^{\mathrm{sh}}$. %
  Further, again by Proposition \ref{p:pplimit}, we obtain a lift strict \'etale locally on $S$, as desired.
\end{pf}

\begin{para}
\label{endpf_lacoh2}
  We prove the surjectivity of $(*)$.
  By Lemma \ref{l:flat}, we may assume that $A$ comes from a principally polarized log abelian variety $A_0$ over a log regular base. 
  Since $\coLie$ is a vector bundle, we reduce to the surjectivity of $(*)$ for $A_0$.
  Hence we may assume that the base $S$ is log regular in $(*)$. 
  Let $P$ be any complete model of $A$. 
  Let $V$ be any fs log scheme which is log \'etale over $A$. 
  Then $V \times_{A}P \to V$ is a log blow-up.
  Since $S$ is log regular,  $V$ is also log regular and we have $R(V \times_{A}P \to V)_*\Ga=\Ga$ 
 (cf.\ \cite{KKMS} p.44, Chapter I, Section 3, Corollary 1 to Theorem 12). 
  From this, we see that $H^m(A,\omega^1_{A/S})=H^m(P,\omega^1_{P/S})$ for every $m$. 
  In particular, $H^0(A,\omega^1_{A/S})=H^0(P,\omega^1_{P/S})$.

  Now take an element $\alpha \in \coLie(A)(S)$. 
  We use the same notation as in \ref{inj}.
  By (1) in \ref{inj}, the image of $\alpha$ in $\coLie(A)(\hat S)$ is in the image of $H^0(A_{\hat S}, \omega^1_{A_{\hat S}/\hat S})$.
  Take a complete model $P$. 
  Since the image of $\alpha$ in $\coLie(A)(\hat S)$ can be lifted to
$H^0(P_{\hat S}, \omega^1_{P_{\hat S}/\hat S})$, by Artin's approximation theorem, 
the image of $\alpha$ in $\coLie(A)(S^{\mathrm{sh}})$ can be lifted to
$H^0(P_{S^{\mathrm{sh}}}, \omega^1_{P_{S^{\mathrm{sh}}}/S^{\mathrm{sh}}})$.
  Then, $\alpha$ itself can be lifted to 
$H^0(P, \omega^1_{P/S})$. 
  But we have $H^0(P, \omega^1_{P/S})=H^0(A, \omega^1_{A/S})$, which implies that $\alpha$ is in the image of $H^0(A, \omega^1_{A/S})$.

  The proof of Proposition \ref{p:lacoh2} is completed. 
\end{para}

\begin{para}
  Let $S$ and $A$ be as in Proposition \ref{p:lacoh2}. 
  Let $\omega^1 = \omega^1_{A/S}$. 

  Since $\cH^0(A,\Gmlog)=\Gmlog$ by \cite{KKN5} Proposition 2.1, we have an exact sequence 
$$0 \to \cH^0(A,\omega^1) \to \cH^1(A, [\Gmlog \to \omega^1]) \to \cH^1(A,\Gmlog) \to \cH^1(A,\omega^1). $$ 
  Here $[\Gmlog \to \omega^1]$ means the complex where $\Gmlog$ is in degree zero and the map is $d\log$. 

  From this, we have an exact sequence %
$$0 \to \cH^0(A,\omega^1)' \to \cH^1(A, [\Gmlog \to \omega^1])' \to \cH^1(A,\Gmlog)' \to \cH^1(A,\omega^1)'.$$
  (See \ref{Extnatural} for the notation $(-)'$.)

  Since $A$ as a group sheaf is locally of finite presentation (\cite{KKN5} Proposition 12.7), there is a natural map 
$A \to \cH^1(A,\Gmlog)'$. 

\end{para}

\begin{prop}\label{ue2} 
  Let $S$ and $A$ be as in Proposition $\ref{p:lacoh2}$. 
The universal additive extension is obtained also as the fiber product 
$E'(A)$ of 
$$\cH^1(A, [\Gmlog \to \omega^1])' \to \cH^1(A, \Gmlog)' \leftarrow A.$$
  We have a universal exact sequence 
$$0 \to \cH^0(A,\omega^1)' \to E'(A) \to A \to 0.$$ 
\end{prop}

  We prove a lemma and a proposition we will use in the proof of Proposition \ref{ue2}.   

\begin{lem}
\label{l:cubic}
  Let $A$ be a weak log abelian variety over a noetherian fs log scheme satisfying the condition $1.4.1$ in {\rm $\mathrm{\cite{KKN5}}$}. 
  Then we have $\Ext(A,\Ga) = H^1(A,\Ga)$. 
\end{lem}

\begin{pf}
  This is an analogue of the cubic isomorphism \cite{KKN5} Theorem 2.2 (c), and the general case is reduced to the constant degeneration case as in the same way as in \cite{KKN5} 12.5.
  In the rest of this proof, we assume that $A$ is with constant degeneration. 
  Then, by \cite{KKN5} 6.4.1, we have $H^i(\widetilde A, \Ga)=H^i(B,\Ga)$ (the notation is as usual as in there). 
  Hence we have the spectral sequence $E_2^{ij}=H^i(Y,H^j(B,\Ga)) \Rightarrow H^{i+j}(A,\Ga)$. 
  Since $Y$ acts on $H^j(B,\Ga)$ trivially for every $j$, there is an exact sequence 
$$0 \to \Hom(Y,\Ga) \to H^1(A,\Ga) \to H^1(B,\Ga).$$
  Taking $H^1(B,\Ga) =\Ext(B,\Ga)$ into account, 
by comparing this with the exact sequence (3) in the proof of Proposition \ref{p:coLie}, we have $\Ext(A,\Ga)=H^1(A,\Ga)$, as desired.
\end{pf}

\begin{prop}
\label{p:lacoh}
  Let $A$ and $A'$ be principally polarized log abelian varieties over a noetherian fs log scheme $S$. 
  Then we have the following. 

$(1)$ The restriction of $\cH^1(A,\Ga)$ to $S_{\et}$ is a vector bundle of rank $\dim A$. 

$(2)$ $H^1(A\times A',\Ga)=H^1(A,\Ga)\oplus H^1(A',\Ga)$. 
\end{prop}

\begin{pf}
(1) By Lemma \ref{l:cubic}, %
this is reduced to Proposition \ref{p:coLie}. 

(2) It is by $\cH^0(A,\Ga)=\cH^0(A',\Ga)=\Ga$ (\cite{KKN5} Proposition 12.1) and (1).  %
\end{pf}

\begin{para}
  To have the exact sequence in Proposition \ref{ue2}, %
it is enough to show 
that the composition $A\to \cH^1(A, \Gmlog) \to \cH^1(A, \omega^1)$ is the zero map. 

  The map $A\to \cH^1(A, \Gmlog) \to \cH^1(A, \omega^1)$ %
comes from an element $\alpha\in H^1(A\times A, \mathrm{pr}_1^*\omega^1_{A/S})$ such that (we denote the second $A$ in $A\times A$ by $A'$) it sends a section $s$ of $A'$ to the image of $\alpha$ by the pullback with respect to $s$. 
  By Propositions \ref{p:lacoh} (2) and \ref{p:lacoh2}, we have
\begin{align*}
H^1(A \times A', \mathrm{pr}_1^*\omega^1_{A/S}) 
&= H^1(A\times A', \Ga) \otimes H^0(A\times A', \mathrm{pr}_1^*\omega^1_{A/S})\\ 
&=(H^1(A, \Ga) \oplus H^1(A', \Ga)) \otimes \coLie(A)(S) \\
&=H_1 \oplus H_2, 
\end{align*}
where $$H_1= H^1(A, \Ga) \otimes \coLie(A)(S), \quad H_2=H^1(A', \Ga) \otimes \coLie(A)(S).$$ 
Write $\alpha=\alpha_1+\alpha_2$ with $\alpha_i\in H_i$ $(i=1,2)$. 
  Then the pullback of $\alpha$ in $\cH^1(A, \omega^1)$ under every $s\in A'$ is $\alpha_1$ ($\alpha_2$ is killed by the pullback). That is, the map $A' \to \cH^1(A, \omega^1)$ is a constant map. Because it is a homomorphism, it is the zero map. 
\end{para}

\begin{para}
\label{pf2-ue2}
  We prove Proposition \ref{ue2}. 
  By Proposition \ref{p:lacoh2}, $\cH^0(A,\omega^1_{A/S})'$ is a vector bundle. 
  Hence, by the universality, we have a map $f\colon E(A) \to E'(A)$. 
  To prove that it is bijective, %
it is enough to show that %
the induced map $\coLie(A) \to \cH^0(A,\omega^1_{A/S})'$ by $f$ is the inverse of the map given in the proof of Proposition \ref{p:lacoh2}. 
  By Lemma \ref{l:flat}, we reduce to the case where the base is log regular. 
   Hence we reduce to the case where the log structure of the base is trivial, that is the classical case.  
\end{para}

\begin{para}
  Now we prove Proposition \ref{p:E(A)unique}.
  We construct a map 
$\cExt^{\natural}(A,\Gmlog)' \to E(A')$. 
  Since $A'$ is log smooth, by covering $A'$ by log smooth objects, for each section of 
$\cExt^{\natural}(A,\Gmlog)'$, we can associate a class of $\cH^1(A', [\Gmlog \overset {d\log} \to \omega^1_{A'/S'}])'$. 
  We prove that the image of this class in $\cH^1(A',\Gmlog)'$ belongs to 
$A'{}^*$ so that this class is in 
$E'(A')= E(A')$ (Proposition \ref{ue2}). 
  By definition of $\cExt^{\natural}(A,\Gmlog)$, the image of this class in $\cH^1(A,\Gmlog)$ belongs to $A^*$. 
  Hence we reduce to the following claim. 

\smallskip

\noindent 
{\bf Claim 1.} 
  The map from $A'$ to the fiber product of 
$$\cH^1(A',\Gmlog) \to \cH^1(A,\Gmlog) \leftarrow A$$
is bijective. 

\smallskip

  First we assume that $A'$ is with constant degeneration. 
  By the definition of the dual, it suffices to show that 
the map from $\cExt(A',\Gmlog)$ to the fiber product of 
$$\cH^1(A',\Gmlog) \to \cH^1(A,\Gmlog) \leftarrow \cExt(A,\Gmlog)$$
is bijective. 
  Since $\cH^1(A',\Gmlog)/\cExt(A',\Gmlog)$ injects into $\cHom(A',A'{}^*)$ by \cite{KKN5} Theorem 8.5 (2), it is enough to see that a homomorphism $A' \to A'{}^*$ is zero if the induced homomorphism $A \to A^*$ is zero.  
  This is by Lemma \ref{l:endo}. 

  Next, let $U' \in (\fs/S')$ be of finite type over $S'$ (so of finite type over $\bZ$). 
  We want to prove that an element of $\cH^1(A',\Gmlog)(U')$ belongs to $A'(U')$ if its image in $\cH^1(A,\Gmlog)(U)$ belongs to $A(U)$. 
  By \cite{KKN5} Lemma 12.9, 
we can replace $U'$ by its strict localization. 
  By Artin's approximation theorem, we can further replace $U'$ by the completion. 
  To reduce the problem to what we already see in the constant degeneration case, 
  it is sufficient to use the injective GAGF for $H^1(A',\Gmlog)$ (\cite{KKN6} Theorem 1.2) and 
the GAGF for $A'(S)$ (cf.\ the proof of Proposition \ref{p:ppinj}). 
  We conclude that the element concerned belongs to $A'(U')$, as desired. 

  Thus we have a homomorphism 
$\cExt^{\natural}(A,\Gmlog)' \to E(A')$ which fits into the commutative diagram
$$
\begin{CD}
0 @>>> K @>>> \cExt^{\natural}(A,\Gmlog)' @>>> A'{}^* @>>> 0 \\
@.@VVV@VVV@| \\
0 @>>> \cH^0(A',\omega^1_{A'/S'})' @>>> E(A') @>>> A' @>>> 0
\end{CD}
$$
with exact rows, where $K=\Ker(\cExt^{\natural}(A,\Gmlog)' \to A'{}^*)$.
  Hence it is enough to show that $K \to \cH^0(A',\omega^1_{A'/S'})'$ is bijective.
  By Claim 1, $K$ coincides with $\Ker(\cH^1_{\mathrm{logcrys}}(A,\Gmlog)' \to \cH^1(A',\Gmlog)')$.

\smallskip

\noindent 
{\bf Claim 2.} There is a natural exact sequence 
$$0 \to \cH^0(A',\omega^1_{A'/S'})' \to \cH^1_{\text{logcrys}}(A,\Ga)' \to \Lie(A').$$

\smallskip

  By comparing this exact sequence with the upper row of the above commutative diagram, %
we see that %
the map $K \to \cH^0(A',\omega^1_{A'/S'})'$ 
is identified with the identity, which completes the proof.

  We prove Claim 2.  By \cite{Kato:FI1} Theorem 6.4, we have $\cH^1_{\text{logcrys}}(A,\Ga)=\cH^1(A',\omega^{\cdot}_{A'/S'})$. 
  Further, by Proposition \ref{p:coLie} and Lemma \ref{l:cubic}, we have a natural isomorphism $\Lie(A') \cong \cH^1(A',\omega^0_{A'/S'})'$.
  Hence, to obtain the claimed exact sequence as a part of the log Hodge-de Rham spectral sequence, it is enough to show that the homomorphism 
$\cH^0(A',\omega^1_{A'/S'})' \to \cH^0(A',\omega^2_{A'/S'})$ is zero. 
  Since 
$\cH^0(A',\omega^1_{A'/S'})' \cong \coLie(A')$ (Proposition \ref{p:lacoh2}) is a vector bundle, by Lemma \ref{l:flat}, we may assume that the base is log regular. 
  Then it is zero because it is so over the nonlog open set.
\end{para}

\section{Deformation of log abelian varieties}
\label{s:deformation}
  We discuss the deformation of log abelian varieties, which is used in the proofs of main results. 

\begin{prop}
\label{p:deform}
  Let $S=\Spec R \to S'=\Spec R'$ be a strict nil immersion of affine fs log schemes by a finitely generated and square-zero ideal $I \subset R'$.
  Let $A$ be a principally polarized log abelian variety over $S$. 
  Assume that a lift $A'$ of $A$ over $S'$ is given. 
  Then the following hold. 

$(1)$ The set of the liftings of $A$ over $S'$ is naturally bijective to the underlying set of 

$$\Hom_{\sym}(\coLie(A), \Lie(A))\otimes_{R}I,$$
where $\Hom_{\sym}$ means the subgroup consisting of the self-dual homomorphisms.

$(2)$ Assume further that $R'=R[I]$ and that the lift $A'$ is the canonical one. 
  Then the bijection in $(1)$ is an isomorphism of $R$-modules. 
\end{prop}

\begin{pf}
  (1) Regard $\Lie(A')$ as an $\cO_{S'}$-module on the small \'etale site $\upc S{}'_{\et}$. 
  Tensoring it with the exact sequence $0 \to I \to \cO_{S'} \to \cO_S \to 0$, we have an exact sequence 
$$ 0 \to \Lie(A) \otimes_{\cO_S}I \to \Lie(A') \to \Lie(A) \to 0$$ 
on $\upc S{}'_{\et}$. 
  Let $f\colon \coLie(A) \to \Lie(A) \otimes I$ be a homomorphism of $\cO_{S'}$-modules on $\upc S{}'_{\et}$, which comes from an element of the set concerned. 
  Consider the composite 
  $$\coLie(A') \to \coLie(A) \overset f \to \Lie(A) \otimes I \hookrightarrow \Lie(A').$$ 
  Locally, it lifts to a homomorphism $\tilde f\colon \coLie(A') \to \Lie(E(A'))$ that induces the inclusion $\coLie(A) \to \Lie(E(A))$. 
  Regard $\tilde f$ as a homomorphism of abelian sheaves on the big site $(\fs/S')_{\et}$. 
  Then the image in $E(A')$ of the image of $\tilde f$ is glued to give a subsheaf of $E(A')$ and the quotient $A''$ of $E(A')$ with respect to it does not depend on the choices up to isomorphisms. 

  Further, $A''$ is a log abelian variety on $S'$ which lifts the original $A$. 
  To see it, first we see that $A''$ is log smooth over $S'$. 
  Next, there is a unique surjective homomorphism $A'' \to A'/G'$ which lifts $A \to A/G$.
  Here $G$ and $G'$ are the semiabelian parts of $A$ and $A'$, respectively. 
  Let $G''$ be the kernel of this surjection. 
  Then $G''$ is a strict smooth group object over $S'$ so that it is semiabelian because $G$ is so. 
  Since the pullback of $A''$ to $S$ is $A$, the fibers of $A''$ are log abelian varieties.
  Finally we deduce that $A''$ is separated from the separability of $E(A')$.

  The self-duality implies that $A''$ is principally polarized. 
  Thus we have a map from the set concerned to the set of the liftings, which is injective by the universality of $E(A')$. 
  To prove its surjectivity, it is enough to show that for any lift $A''$, we have $E(A') \cong E(A'')$. 
  To see this, by Proposition \ref{p:pplimit}, we may assume that $S'$ is of finite type over $\bZ$, and we use Proposition \ref{p:E(A)unique}.
  We complete the proof of (1). 

  (2) is a formal consequence from (1). 
\end{pf}

\begin{prop}
\label{p:lift}
  Let $S \to S'$ and $A$ be as in the previous proposition. 
  Then $A$ lifts to $S'$ Zariski locally on $S$. 
\end{prop}

\begin{pf}
  The case where $A$ is with constant degeneration is deduced from the log smoothness of the local moduli in \cite{KKN6} Theorem 2.4. 
  Hence, the case where $S'$ is the spectrum of an Artin ring is proved. 
  Then the case where $S'$ is the spectrum of a complete noetherian local ring is proved by GAGF (\cite{KKN6} Theorem 6.1 and Remark 6.1.1).

  We prove the case where the underlying scheme of $S'$ is a strict localization of an fs log scheme of finite type over $\bZ$. 
  By the previous case, we have a principally polarized log abelian variety over the completion of $S'$. 
  Then by Proposition \ref{p:pplimit} and Artin's approximation theorem, this case follows. 

  Consider the general case. 
  We may assume that $S'$ is of finite type over $\bZ$ by Proposition \ref{p:pplimit}. 
  Again by Proposition \ref{p:pplimit} and the case already proved, we can find a lift \'etale locally on $S$. 
  By Proposition \ref{p:deform}, the obstruction for the existence of a global lift is in $H^1(S, \cHom_{\sym}(\coLie(A),\Lie(A))\otimes I)$, where $I$ is as in Proposition \ref{p:deform}.
  This vanishes Zariski locally. 
  Hence $A$ lifts Zariski locally on $S$. 
\end{pf}

\section{Some related functors}
\label{s:related}

  We introduce some related functors which we will use in the proofs. 

  Throughout this section, let $g$, $W$, $n$ be as in Section \ref{sec:functor}, 
$\Sigma$ an admissible cone
decomposition of $S_{\Q}(W)$ (Definition \ref{defn:adm-decomp}), and 
$$F:= F_{g, n}, \quad F_{\Sigma}:=F_{g,n,\Sigma}\colon (\fs/\bZ[1/n])\to (\mathrm{set})$$ the moduli
functors in Definition \ref{defn:functor}.

  By the next lemma, $F$ is a sheaf with respect to the classical \'etale topology if $n \geq 3$. 

\begin{lem}\label{lem:aut-triv}
  Let $U$ be an fs log scheme, let $n\ge 3$, and let $A$ be a log
  abelian variety over $U$ with a principal polarization $\phi$ and with a
  level $n$ structure $e$. Then $\Aut(A,\phi,e)=\{1\}$.
\end{lem}
\begin{pf}
  By Proposition \ref{p:pplimit}, \cite{KKN4} Proposition 9.2 (1), and the representability of torsion points of a log abelian variety (\cite{KKN4} Proposition 18.1 (1)), we can spread out a triple $(A,\phi,e)$ and an automorphism of it. 
  Hence we may assume that $\upc U$ is of finite type over $\bZ$. 
  Since $A$ is locally of finite presentation (\cite{KKN5} Proposition 12.7), we may replace $U$ by a strict localization. 
  By Artin's approximation theorem and again by \cite{KKN4} Proposition 9.2 (1), we may further replace $U$ by the spectrum of a complete noetherian local ring. 
  Then by GAGF for a log abelian variety (\cite{KKN6} Theorem 6.1), we may replace $U$ by the spectrum of an Artin local ring. 
  Lastly, by Lemma \ref{l:endo}, an automorphism of $A$ whose restriction to the closed point is the identity is the identity so that we may replace $U$ by the closed point. 

  The rest is to prove the case where $\upc U$ is the spectrum of a field. 
  First we prove that $\Aut(A, \phi, e)$ is torsion-free. 
  Since $n\ge 3$, as in the classical case, it is enough 
to show that for a prime number $\ell$ which is invertible on $U$ and two log abelian 
varieties $A_1$ and $A_2$ over $U$, the natural map 
$\Hom(A_1, A_2) \to \Hom_{\Z_\ell}(T_\ell A_1, T_\ell A_2)$ is 
injective.
  Let $f \in \Hom(A_1, A_2).$ 
  Then the induced $T_\ell(f)$ preserves the weight filtrations 
(which means the filtration defined by the homomorphisms in \cite{KKN4} 18.9.1 and 18.9.2).
  Thus the above is reduced to the classical case 
(cf.\ \cite{Mumford:textbook} p.176, Theorem 3) and hence $\Aut(A,\phi, e)$ 
is torsion-free.
  Note here that $\Hom(A_1, A_2)$ is also finitely generated. 
  This is also reduced to the classical case as follows. 
  Since $A_i$ is with constant degeneration, we can consider the associated log 1-motif 
$[Y_i \to G_{i,\log}]$ to $A_i$ ($i=1, 2$). 
  Then we have 
$\Hom(A_1, A_2) \hookrightarrow \Hom(Y_1, Y_2) \times 
\Hom(G_{1,\log}, G_{2,\log}) =
\Hom(Y_1, Y_2) \times \Hom(G_1, G_2)$, 
where the equality is by \cite{KKN2} Proposition 2.5.
  Hence, $\Hom(A_1, A_2)$ is finitely generated.
  In particular, $\End(A)$ is finitely generated.

  Next, we prove that $\Aut(A,\phi)$ is finite, which competes the proof.
  We may assume that $Y$ for $A$ is constant. 
  We define the Rosati involution as in the classical case: 
for $g \in \End(A)$, let ${}^*g$ be the homomorphism 
$A\overset {\phi} \to A^* \overset {g^*} \to A^* 
\overset {{\phi}^{-1}} \to A$. 
  Then ${}^*(\cdot)$ gives an anti-automorphism of $\End(A)$. 
  Let $\langle g_1, g_2\rangle$ be the trace of the induced 
endomorphism on the Tate module by $g_1\compo{}^*g_2$. 
  Since this trace is the sum of the traces of the 
endomorphisms on $Y$, $X$ and $B$, where $X \times Y \to 
M_S^{\gp}/{\cal O}_S^{\times}$ is the associated pairing and 
$B$ is the abelian part, $\langle g_1, g_2\rangle$ defines 
a bilinear form on $\R\otimes\End(A)$ as in the classical 
case.
  Further, by the definition, $\langle g, g \rangle =1$ 
for any automorphism $g$. 
  The rest is to show the positivity, that is, $\langle g,g\rangle = \mathrm{Tr}
(g\compo{}^*g)>0$ for $g\not=0$. 
  As is explained above, the trace is the sum of the three 
traces.
  The classical part is known (cf.\ \cite{Mumford:textbook} 
p.192, Theorem 1). 
  The part on $Y$ and the part on $X$ can be treated 
similarly. 
  We write down the former case.
  Let $f\colon Y_{\Q} \to Y_{\Q}$ be the induced 
homomorphism from $g$ and ${}^*f\colon Y_{\Q} \to 
Y_{\Q}$ induced from ${}^*g$. 
  We have to prove $\mathrm{Tr}(f\compo{}^*f)>0$ if $f\not=0$. 
  Let $\phi\colon Y_{\Q} \to X_{\Q}$ be the homomorphism 
induced by $\phi$. 
  Then ${}^*f$ coincides with $\phi^{-1}f\spcheck \phi$, 
where $f\spcheck\colon X_{\Q} \to X_{\Q}$ is the dual of $f$. 
  Let $\langle y, z\rangle_Y$ be $\langle \phi(y), z\rangle 
\in M_S^{\gp}/{\cal O}_S^{\times}$.
  Then 
$\langle y, f(z) \rangle_Y=
\langle \phi(y), f(z) \rangle=
\langle f\spcheck\phi(y), z\rangle=
\langle \phi^{-1}f\spcheck\phi(y), z\rangle_Y=
\langle {}^*f(y), z\rangle_Y$. 
  Let $N\colon M_S^{\gp}/{\cal O}_S^{\times}\to \N$ 
be a homomorphism whose kernel is trivial. 
  Then $N(\langle\cdot, \cdot \rangle_Y)$ is a positive definite form. 
  Since $f$ and ${}^*f$ are adjoint to each other with 
respect to this form, $\mathrm{Tr}(f\compo{}^*f)$ is positive. 
\end{pf}

\begin{para}\label{para:I/Gamma}
  To prove main theorems, we define 
the sheaf $F'$ on $(\fs/\Z[1/n])_{\et}$ as 
\begin{align*}
F'(U):=\{&
(A, \phi, e, f)\}/\cong, 
\end{align*}
$U \in (\fs/\Z[1/n]),$ 
where $(A,\phi,e)$ is a $g$-dimensional principally polarized log abelian variety over $U$ with a level $n$ structure and 
$f$ is a surjective homomorphism $W \to \ol Y(A)$.  

  This $F'$ is indeed a sheaf because $F$ is so.  
\end{para}

\begin{prop}\label{sbprop:F'toF}
  The forgetful morphism $F' \to F\ ; \ 
(A, \phi, e, f) \mapsto (A, \phi, e)$ is represented by 
strict \'etale surjective algebraic spaces.  
  More strongly, it is represented by surjective Zariski local isomorphisms. 
\end{prop}

\begin{pf}
  Let $(A, \phi, e)\colon U \to F$ be a morphism from any 
fs log scheme $U$ over $\bZ[1/n]$. 
  The functor $F' \times_FU$ over $U$ is represented by the sheaf 
of surjective homomorphisms $W \to \ol Y(A)$. 
  Let $u \in U$, put $Y := \ol Y(A)_{\bar u}$, and, shrinking $U$ if necessary, we may assume that there is a surjection $Y \to \ol Y(A)$.
  Then $F' \times_FU$ is represented by the fs log scheme 
$(\coprod_h U_h)/\sim$, where $h$ runs over the set of the 
surjections $W \to Y$, $U_h$ is a copy of $U$, and $U_h$ 
and $U_{h'}$ are glued by the open subscheme on which 
$h$ and $h'$ induce the same $W \to \ol Y(A)$. 
  The formation is compatible with any base change so that 
$F' \to F$ is globally represented. 
\end{pf}

\noindent {\it Remark.}
  Further, $F' \to F$ is represented by quasiseparated morphisms. 
  This is seen as follows.
  In the notation in the above proof,  
$\ol Y(A)$ locally comes from an fs log scheme $U$ that is of 
finite type over $\Z$. 
  Hence the glued schemes are quasiseparated. 

\begin{para}
\label{cI}
  Next we introduce the sheaf 
$\cI$ on $(\fs/\Z[1/n])_{\et}$ defined as 
\begin{align*}
\cI(U):=\{&\text{symmetric bilinear form\ }W\times W
\overset b \to (\Gmlog/\Gm)(U) \\
&\text{\ 
which is positive semi-definite}\},\quad 
\end{align*}
$U \in (\fs/\Z[1/n]),$ 
where $b$ is said to be {\it positive semi-definite} if, for 
any $u \in U$ and any homomorphism $(M_{U}/\cO_U^{\times})_{\bar u} \to \R^{(+)}_{\ge 0}$, the induced $W \times 
W \to \R$ is positive semi-definite. 
  Here $\R^{(+)}_{\ge 0}$ means the set $\{a \in \R\,\vert\,
  a\ge 0\}$ regarded as a monoid with respect to the addition.

  We define the canonical morphism 
$$F' \to \cI$$ 
of sheaves by associating 
to $(A, \phi, e, f)$ over $U \in (\fs / \Z[1/n])$ 
a bilinear form 
$$
W\times W \overset {f \times f} \to \ol Y(A) \times 
\ol Y(A)  \overset {\phi \times \mathrm{id}} \to
\ol X(A) \times \ol Y(A)  
\overset {\langle \ , \ \rangle}\to \Gmlog /\Gm.$$

  We construct an explicit covering of $\cI$ by representable objects. 
  Let $\cT$ be the disjoint union of 
$$V_{\sigma} := \Spec \bigl((\bZ[\tfrac 1n])[\sigma{}^{\vee} \cap \mathrm{Sym}^2_{\Z}
(W)]\bigr),$$ where $\sigma$ runs over the set $K$ of all sharp 
finitely generated $\Q$-cones consisting of positive 
semi-definite symmetric bilinear forms $W \times W \to \bR$.

  We have the natural morphism 
$$\cT \to \cI.$$ 
\end{para}

\begin{prop}\label{sbprop:I}
  The natural morphism 
$\cT \to \cI$ is surjective as a morphism of sheaves, 
$\cT\times_{\cI}\cT$ is representable, the 
two projections $\cT\times_{\cI}\cT \rightrightarrows \cT$ are log \'etale and the 
natural morphism $\cT\times_{\cI}\cT \to \cT\times_{\bZ[\tfrac {\scriptscriptstyle 1}{\scriptscriptstyle n}]}\cT$ is 
of finite type. 
\end{prop}

\begin{pf}
  We prove that $\cT \to \cI$ is surjective. 
  Let $b$ be an element of $\cI(U)$ $(U \in (\fs/\bZ[1/n]))$ and we prove that $b$ comes from a section of $\cT$. 
  We may assume that $U$ has a chart. 
  Let $\sigma$ be the dual of the inverse image of $M_U/\cO_U^{\times}$ by the homomorphism 
$\mathrm{Sym}^2_{\Z}(W) \to \Gmlog/\Gm$ induced by $b$. 
  Then, strict \'etale locally on $U$, $b$ comes from a section of $V_{\sig}$. 
  Hence, $\cT \to \cI$ is surjective. 

  Let $\sigma$, $\tau \in K$. 
  To prove the rest, it is enough to show that 
$R:= V_{\sigma} \times_{\cI}V_{\tau}$ is representable, that 
$R \to V_{\sigma}$ is log \'etale and that $R$ is 
quasicompact over $\bZ[\tfrac 1n]$. 
  But the map 
$V_{\sigma\cap \tau}\times \mathbb G_m 
\to R$ sending $(t, u)$ to $(t, t\cdot u)$, 
where we denote by the same symbol the images of $t$, 
is an isomorphism, which suffices.  
\end{pf}

  Let $F'_\cT = F'\times_\cI\cT$. 
  The next proposition is an essence of the proofs of our main results.
  Once it is proved, all results are easily deduced from it. 

\begin{prop}\label{prop:classicalArtin}
  The functor 
$F'_\cT$, restricted 
on $(\mathrm{sch}/\oc \cT) \subset (\fs/\cT)$, is represented by 
a separated algebraic space over $\oc \cT$. 
\end{prop}

  We prove this in the next section. 

\section{Artin's criterion}
\label{s:Artin}
  In this section, we prove Proposition \ref{prop:classicalArtin} by the classical Artin criterion. 
  This section is a core of the proofs in this part of this series of papers. 

\begin{lem}\label{sblem:fppf}
  The 
quotient sheaf $\Gmlog/\Gm$ is an fppf sheaf. 
\end{lem}

\begin{pf}
  This is proved in \cite{Kato:FI2} 3.5. 
  We include the proof for readers' conveniences. 
  Let $U' \to U$ be a strict fppf covering. 
  Let $U'':= U'\times_UU'$ and it is enough to show that 
$$\Gamma(U, \Gmlog/\Gm) \to \Gamma(U', \Gmlog/\Gm) \rightrightarrows 
\Gamma(U'', \Gmlog/\Gm)$$ is exact. 
  We may assume that $U$ has a chart. 
  Then every $\Gmlog/\Gm$ is the inverse image of $\Gmlog/\Gm$ on $U_{\mathrm{Zar}}$. 
  Hence the above exactness is reduced to the exactness of the diagram of topological 
spaces $U'' \rightrightarrows U' \to U$. 
\end{pf}

\begin{para}
\label{formal_smoothness}
  We prove Proposition \ref{prop:classicalArtin}.
  The functor concerned over $(\mathrm{sch}/\oc \cT)$ associates to 
each $\oc \cT$-scheme $U$ the set of isomorphism classes of 
$(A, \phi, e, f)$ over $(U, M_\cT)$ whose induced 
section in 
$\cI((U, M_\cT))$ 
coincides 
with the pullback of the one 
determined by the given $\cT \to \cI$. 
  Here $M_\cT$ is the inverse 
image log structure of $\cT$.
  We denote this functor by $H$. 

  By the classical Artin criterion \cite{Artin:criterion} Theorem 5.3, to prove Proposition \ref{prop:classicalArtin}, it is enough to show that the eleven conditions $[0']$--$[5']$(c) there are satisfied for $H$. 

  Henceforth in this proof, the notation follows loc.\ cit.
  As a deformation theory (loc.\ cit.\ Definition 5.2) for $H$, 
we take the standard one : $D (A_0, M, \xi_0) = 
H_{\xi_0}(A_0 [M])$. 
  To see that it is actually a deformation theory, it is enough to 
check the bijectivity of the map $(*)$ in loc.\ cit.\ 
p.48, which also implies the conditions [$4'$](b) and [$5'$](a) as 
explained there. 
  This bijectivity is proved by 
the calculation of deformations of log abelian varieties as 
follows. 
  Let us recall the situation.  First let $A_0$ be a noetherian 
$\O_\cT$-integral domain. 
  Let $A' \to A \to A_0$ be infinitesimal extensions. 
  Let $B \to A_0$ be another infinitesimal extension.  
  For every map $B \to A$, consider the map 
\begin{equation}\tag{$*$}
H(A'\times_AB) \to H(A') \times_{H(A)} H(B).
\end{equation}
  We must show that it is bijective. 
  Let $(A_1, \phi_1, e_1, f_1) \in H(B)$. 
  By Proposition \ref{p:deform} (1) and Proposition \ref{p:lift}, we have a bijection from the inverse 
image of $(A_1, \phi_1, e_1, f_1)$ by $H(A' \times_AB) \to 
H(B)$ to 
$$\Hom_{\sym}(\coLie(A_1), \Lie(A_1)) \otimes I,$$
where $I= \Ker(A'\times_AB \to B)$.
  Since 
$I=\Ker (A' \to A),$ 
again by 
Proposition \ref{p:deform} (1) and Proposition \ref{p:lift}, it is naturally bijective to the inverse image of the image in $H(A)$ of 
$(A_1, \phi_1, e_1, f_1)$ by $H(A') \to H(A)$. 
  Hence $(*)$ is bijective. 
\end{para}

\begin{para}
   We prove [$0'$], that is, that $H$ is a sheaf for the 
fppf topology. 

  It is enough to show that $F'_\cT$ is an fppf sheaf.  
  Let $p\colon U' \to U$ be a strict fppf covering of a 
$\cT$-fs log scheme.  
  We prove that $F'_\cT(U) \to F'_\cT(U')$ is injective. 
  Let $(A_1, \phi_1, e_1, f_1)$, $(A_2, \phi_2, e_2, f_2) \in F'_\cT(U)$ and 
assume that $p^{-1}(A_1, \phi_1, e_1, f_1) \cong 
p^{-1}(A_2, \phi_2, e_2, f_2)$. 
  We will prove $(A_1, \phi_1, e_1, f_1) \cong (A_2, \phi_2, e_2, f_2)$. 
  Since $F'_\cT$ is an \'etale  %
sheaf, we may assume that there is 
an admissible pairing $b \colon W \times W \to \mathcal S^{\gp}$ 
$(\mathcal S$ is an fs monoid) with a homomorphism ${\mathcal S} \to 
(\Gmlog/\Gm)_U$ such that the composite $W \times W \to (\Gmlog/\Gm)_U$ 
is the given pairing $U \to \cI$ 
which induces $A_1/G_1 \cong A_2/G_2 \cong \HOM(X, \Gmlog/\Gm)^{(Y)}/
\ol Y$. 
  Here $X = Y = W$, and $G_1$ and $G_2$ are the semiabelian parts 
of $A_1$ and $A_2$ respectively.
  The assumption together with Lemma \ref{lem:aut-triv} implies 
that there is an isomorphism $p^{-1} A_1
\cong p^{-1} A_2$ whose two pullbacks by the 
projections $U'':=U'\times_UU' \rightrightarrows U'$ coincide. 
  It induces an isomorphism $p^{-1} \widetilde A_1
\cong p^{-1} \widetilde A_2$ whose two pullbacks by the 
projections $U''\rightrightarrows U'$ coincide. 
  Now let $C \subset \HOM({\mathcal S}, \N) \times 
\HOM(X, \Z)$ be as in \cite{KKN4} 2.2.
  Let $C(m)$ be the cones such that $C = \bigcup C(m)$ as in \cite{KKN1} 3.4.9.   By taking $C(m)$-models, we have an isomorphism 
$p^{-1} \widetilde A_1^{(C(m))}\cong p^{-1} \widetilde A_2^{(C(m))}$ whose two pullbacks by the 
projections $U''\rightrightarrows U'$ coincide. 
  Then by the fppf descent of algebraic spaces 
together with the descent of the fs log structure 
(cf.\ \cite{Olsson} Appendix Corollary A.5), we see that the 
isomorphism of models descends into an isomorphism $\widetilde A_1^{(C(m))}
\cong \widetilde A_2^{(C(m))}$ of algebraic spaces with fs log structures. 
  By taking the union, we have $\widetilde A_1\cong \widetilde A_2$. 
  This isomorphism preserves the group structure since it is so fppf locally. 
  Further it induces $\overline Y_1 \cong \overline Y_2$ since it is so fppf locally, where $\overline Y_i = \Ker(\widetilde A_i \to A_i)$ ($i=1,2$). 
  Thus we have $A_1\cong A_2$. 
  This isomorphism respects polarizations, level structures and surjections. 
  This is by the fact that $A_i$ $(i=1, 2)$ is an fppf sheaf, which is seen by observing 
that $\tilde A_i$ is a k\'et sheaf (\cite{KKN4} Theorem 5.1) and 
that 
$\tilde A_i$ is the union with respect to the k\'et topology of log algebraic spaces in the first sense.
\end{para}

\begin{para}
  Next we prove the exactness of $F'_\cT(U) \to F'_\cT(U') \rightrightarrows F'_\cT(U'')$.  
  Let $(A', \phi', e', f')$ be in the difference kernel of 
$F'_\cT(U') \rightrightarrows F'_\cT(U'')$.  
  We will prove that $(A', \phi', e', f')$ comes from $F'_\cT(U)$. 
  Note that the cocycle condition holds by Lemma \ref{lem:aut-triv}. 
  Again by that $F'_\cT$ is an \'etale %
sheaf and that we already showed that 
$F'_\cT$ is a separated presheaf for 
the fppf topology, we may assume that 
there is a pairing $b \colon W \times W 
\to \mathcal S^{\gp}$ as before which induces $A'/G' 
\cong \HOM(X, \Gmlog/\Gm)^{(\Sigma)}/\ol Y$. 
  (Here and hereafter the notation is similar to that in the previous 
paragraph.) 
  Again by taking $\tilde A{}'$ and its models, we have descent data of log algebraic spaces in the first sense.
  By the fppf descent of log algebraic spaces in the first sense 
as before, $\widetilde A{}'{}^{(C(m))}$ descends as 
a log algebraic space in the first sense.
  Then $\widetilde A$ descends. 
  Next, the pairing $W \times W \to \Gmlog/\Gm$ descends by Lemma \ref{sblem:fppf}. 
  The descended pairing is admissible since it is so fppf locally. 
  Since the pairing descends, $\overline X$ and $\overline Y$ also descend. 
  The homomorphism $\overline Y \to \widetilde A$ descends.  
  It is injective since it is so fppf locally. 
  Let $A$ be the cokernel of $\overline Y \to \widetilde A$. 

  We prove that $A$ is a log abelian variety. 
  Since $\HOM(X, \Gmlog/\Gm)^{(Y)}/\ol Y$ is an fppf sheaf 
by Lemma \ref{sblem:fppf}, the surjective morphism $A'\to 
\HOM(X, \Gmlog/\Gm)^{(\Sigma)}/\ol Y$ also descends. 
  Hence, the semiabelian part also descends.  
  Thus the exact sequence $$0 \to G \to A \to \HOM(X, \Gmlog/\Gm)^{(Y)}/\ol Y \to 0$$ exists. 
  The diagonal is finite because it is so fppf locally. 

  The principal polarization descends by Proposition \ref{p:ppfppf}, which implies the pointwise polarizability. 
  Once the polarization descends, the level structure and the surjection descend.
  This completes the proof of [$0'$].  
\end{para}

\begin{rem}
  There are two alternative proofs for [$0'$]. 
  One is by the use of \cite{KKN6} Theorem 4.7, that is, the descent of $A$ is reduced to that of models. 
  Another proof is as follows. 
  We can prove [$1'$] and [$2'$] in advance. 
  See \ref{1'} and \ref{2'} below.
  Then, by [$1'$], we may assume that $U$ is of finite type over $\bZ$, and replace $U$ by a strict localization. 
  By Artin's approximation theorem and [$1'$] again, we further replace $U$ by the spectrum of a complete noetherian local ring. 
  Then, by [$2'$], we reduce to the case where $U$ is the spectrum of an Artin local ring so that it suffices to show the fppf descent of log 1-motives $[Y \to G_{\log}]$. 
  The lattice $Y$ descends. 
  By \cite{KKN2} Proposition 2.5, the descent data on $G_{\log}$ give those on $G$. 
  Hence $G$ descends and $G_{\log}$ also descends. 
  The rest is to see that the homomorphism $Y \to G_{\log}$ descends. 
  For it, it suffices to show that $G_{\log}$ is an fppf sheaf. 
  By using the exact sequence $0 \to T_{\log} \to G_{\log} \to B \to 0$, where $T$ and $B$ are the torus and the abelian part, respectively, we reduce the problem to the vanishing of $R^1\eta_*\Gmlog$, where $\eta$ is the projection of the fppf site to the \'etale site. 
  The last vanishing is a part of the log Hilbert 90 (see \cite{Kato:FI2} Corollary 5.2). 
\end{rem}

\begin{para}
\label{1'}
  We prove $[1']$ that $H$ is locally of finite presentation. 
  It is enough to show that $F'$ is locally of finite presentation, that is, 
that $\varinjlim F'(U_i) \to F'(U)$ is bijective, where 
$(U_i)_i$ is an inverse system of affine $\bZ[1/n]$-fs log schemes with 
strict transition morphisms and $U := \varprojlim U_i$. 

  This is by Proposition \ref{p:pplimit} by taking care of level structure and the surjection $W \to \overline Y(A)$. 
  The level structure can be spread out by the representability of the torsion points (\cite{KKN4} Proposition 18.1). 
  It is obvious that the surjection is spread out.
\end{para}

\begin{para}
\label{2'}
  We prove $[2']$. 
  It is enough to show that $F'(R) \to \varprojlim 
F'(R/m^n)$ is injective and has a dense image, where 
$(R, m, k)$ is a complete noetherian local $\O_\cT$-algebra.  
  In fact the above map is bijective. 
  This is essentially by \cite{KKN6} Theorem 6.1. 
  We have to take care of the polarization, the level structure, and the surjection.
  For the polarization, see \cite{KKN6} Remark 6.1.1. 
  The level structure is algebraized because the group of torsion points of a log abelian variety is represented by finite log flat group schemes (\cite{KKN4} Proposition 18.1).
  Finally, the formal surjection $W \to \overline Y(A)$ determines the surjection on $R$. 
\end{para}

\begin{para}
  We prove the condition [$3'$](a). 
  It is enough to show the following: 
  Let $S$ be a geometric discrete valuation ring with an fs log structure. 
  Then $H(S) \to H(K)$ is injective, where $K$ is the fraction field of $\cO_S$. 

  If $S$ has the direct image log structure from $K$, then $H(S) \to H(K)$ is bijective 
because ${\cal C}_0^{\mathrm{ptpol}}={\cal C}_0^{\mathrm{pol}}={\cal C}_2^{\mathrm{pol}}$ (in the notation there) by \cite{KKN6} Theorem 3.4 (1) and (3). 
  (Though $S$ is not necessarily an fs log scheme, we define $H(S)$ here as in \cite{KKN6} 3.1.)

  The general case is reduced to this case as follows. 
  It suffices to prove that $H(S) \to H(S')$ is injective, where $S'$ is the one endowed with the direct image log structure from $K$. 
  Since $H(S) \to H(S^{\rm {sh}})$ is injective by [$0'$] and [$1'$], where $S^{\rm{sh}}$ is the strict henselization of $S$, we can replace $S$ and $S'$ by their strict localizations.  
  Furthermore, since $H(S) \to H(\hat S)$ is injective by [$1'$] and Artin's approximation theorem, where $\hat S$ is the completion of $S$, we can further replace $S$ and $S'$ by their completions.  

  Thus we may assume that $S$ and $S'$ are complete. 
  Then by [$2'$], we can replace $S$ and $S'$ by the spectrum of Artin rings. 
  Hence the log abelian varieties concerned are with constant degeneration. 
  Let $\xi$ and $\eta$ be elements of $H(S)$.
  By definition of $H$, the admissible pairings $X \times Y \to \Gmlog/\Gm$ of $\xi$ and $\eta$ are common. 
  Hence the torus part $T$ is common. 
  Assume that the pullbacks to $S'$ of $\xi$ and $\eta$ coincide. 
  Since the underlying schemes of $S$ and $S'$ are the same, the abelian part $B$ and the semiabelian part $G$ are also common.
  Hence it is enough to show that the homomorphism $Y \to G_{\log}$ is determined by the pullback to $S'$. 
  But the induced homomorphism $Y \to G_{\log}/G = T_{\log}/T$ is common.
  Since there is an exact sequence $0 \to \Hom(Y,G) \to \Hom(Y,G_{\log}) \to\Hom(Y,G_{\log}/G)$ and an element of $\Hom(Y,G)$ is determined by the pullbacks, 
we have $\xi=\eta$. 
\end{para}

\begin{para}
\label{4'a}
  We prove [$4'$](a), which means $H_{\xi_0}(A_0[M])$ commutes with localization in $A_0$ and is a finite module when $M$ is free of rank one. 
  Let $A$ be the polarized log abelian variety corresponding to $\xi_0$. 
  Then the condition [$4'$](a) is reduced to the equality 
$$H_{\xi_0}(A_0[M])=\Hom_{\sym}(\coLie(A), \Lie(A))\otimes_{A_0}M$$
as $A_0$-modules. 

  If we neglect the additional structures, this is by Proposition \ref{p:deform} (2). 
  So it is enough to show that the additional structures deform uniquely.
  By \cite{KKN4} Proposition 18.1 (3), 
the torsion points are k\'et locally isomorphic to $(\Z/n\Z)^{2g}$. 
  Since the k\'et site over $\Spec(A_0)$ and that over $\Spec(A_0[M])$ are equivalent, the level structure uniquely deforms. 
  The surjection also uniquely extends, which completes the proof of [$4'$](a).
\end{para}

\begin{para}
\label{local}
  To see that the conditions [$3'$](b), [$4'$](c) and 
[$5'$](c) are valid for $H$, we may assume 
by \cite{KKN2} Theorem 4.6 (2) that 
the $\xi$, $\eta$ etc.\ concerned in these conditions are with 
constant degeneration because $M_{A_0}/
\O^{\times}_{A_0}$ is locally constant around the generic 
point.
  (Here we use [$4'$](a) for [$4'$](c).)  
  Then their lattices $Y$ are constant because the monodromy acts on $Y$ trivially by Lemma \ref{lem:aut-triv}.
  Let $r$ be their rank. 
  Fix an isomorphism $\Z^r =:Y \cong \ol Y(A')$ 
and a surjection $W \to Y$ 
which 
are compatible with 
$f \colon W \to \ol Y$, and denote by $b$ the map 
$U := (\Spec (A_0), M_\cT) 
\to \cI$. 
  Since $U$ is quasicompact, we can take a finitely 
generated cone 
$\sig$ as in \cite{KKN6} Section 2 such that 
for any $u \in U$ and any homomorphism $M_{U, \bar u}
\to \Q_{\ge0}^{(+)}$, the induced pairing $Y \times 
Y \to \Q$ by $b$ belongs $\sig$. 
  Then $b$ factors through $\Spec(\sig{}^{\vee})$ and further 
through $\cT_0:= \cT \times_{\cI}\Spec(\sig{}^{\vee})$ also after 
localizing $U$. 
  Hence we may assume that $b$ factors through $\cT_0$. 
  (Here we use again [$4'$](a) for [$4'$](c).)  
  Now $\xi$ etc.\ can be lifted to sections of 
$(F_{g,r,n,\sigma})_{\cT_0}$, and we can replace $H$ by 
the restriction $H_{g,r,n,\sigma}$ of $(F_{g,r,n,\sigma})_{\cT_0}$ 
in (sch$/\oc {\cT_0})$. 
  (Note that for [$4'$](c) and [$5'$](c), we use here 
the fact that the deformations of log abelian varieties 
with the constant $\ol Y$ must be with the 
constant $\ol Y$.) 
  Since $(F_{g,r,n,\sigma})_{\cT_0}$ is pro-representable 
by \cite{KKN6} Theorem 2.4, $H_{g,r,n,\sigma}$ 
is also pro-representable (by the underlying formal 
scheme). 
  From this, we can see that [$3'$](b), [$4'$](c) and 
[$5'$](c) are valid for $H_{g,r,n,\sigma}$. 
\end{para}

\begin{para}
  Since we already see in \ref{formal_smoothness} that the conditions [$4'$](b) and [$5'$](a) are satisfied, the rest is [$5'$](b). 
  Since $A_0 \times_K A'_K \to A_0\times_K A_K$ in loc.\ cit.\ (the notation is as in there) is a 
nilpotent thickening of local rings with square zero ideal, 
this is reduced to Proposition \ref{p:lift} as follows. 
  In fact, we prove that, in general, for a nil immersion $\Spec R \to \Spec R'$ of affine local schemes over $\oc \cT$ defined by a (not necessarily finitely generated) 
ideal $I$ with $I^2= 0$, any element of $H(R)$ can be lifted to $H(R')$. 
  Take an inductive system $(R'_i)_i$ of noetherian local sub $\cO(\oc \cT)$-algebras of $R'$ such that $\varinjlim R'_i$ is isomorphic to $R'$. 
  For any $i$, let $R_i$ be the image of $R'_i \to R$ so that $\varinjlim R_i=R$.
  By [$1'$], any element of $H(R)$ can be lifted to an element of $H(R_i)$ for some $i$.
  Hence we may assume that $R'$ is noetherian. 
  In this case, by Proposition \ref{p:lift}, any polarized log abelian variety over $R'$ can be lifted to $R$. 
  The level structure and the surjection uniquely extend (cf.\ \ref{4'a}), which completes the proof of [$5'$](b).

  Thus the proof of Proposition \ref{prop:classicalArtin} is completed. 
\qed
\end{para}

\section{Representability of moduli functors}
\label{sec:pf-main}
  In this section, we return to the study of the moduli functors of log abelian varieties defined in Section \ref{sec:functor}, and prove Theorems \ref{thm:mainthm1} and \ref{thm:mainthm2} except properness based on Proposition \ref{prop:classicalArtin}.
\begin{para}
\label{logreg}
  We deduce from Proposition \ref{prop:classicalArtin} we have shown in the previous section, that $F'_\cT$ is represented by a log smooth log algebraic space over $\cT$ in the first sense. 

  In fact, Proposition \ref{prop:classicalArtin} gives us an underlying algebraic space, denoted by $\oc X$, over $\oc \cT$. 
  It contains as a dense open subset the part which is the moduli of abelian varieties without degeneration. 
  We denote this open set by $U$. 
  Endow $\oc X$ with the inverse image log structure from $\cT$. 
  Note that by the local theory (\cite{KKN6} Section 2), it coincides with the log structure by the complement of $U$ and the resulting log algebraic space $X$ in the first sense is log smooth over the base $\bZ[1/n]$.

  By definition, there is a canonical morphism $X \to F'_\cT$.
  It suffices to show that it is an isomorphism.
  Since each log locus (i.e., a subscheme where $M/\cO^{\times}$ is locally constant) of $X$ is a moduli of log abelian varieties with constant degeneration (\cite{KKN6} Section 2, cf.\ \ref{local}), the canonical morphism is formally an isomorphism. 
  Hence, by the argument in \cite{KKN3} 5.4, we reduce to the fact that $F'_\cT$ has the two properties that $F'_\cT$ is locally of finite presentation and that $F'_\cT$ has the uniqueness for GAGF. 
  These are reduced to that $F'$ has the same two properties, which is shown in \ref{1'} and in \ref{2'}, respectively.
\end{para}

\begin{para}
\label{F'_sig}
  {\sc Proof of Theorem \ref{thm:mainthm2} except properness}. 
  We prove a slight variant of Theorem \ref{thm:mainthm2} as follows, which implies Theorem \ref{thm:mainthm2} except properness.
  Let $\Sig$ be a fan in $S_{\bQ}(W)$ supported by positive semi-definite bilinear forms which is not necessarily complete and not necessarily stable under the action of $\Aut_{\bZ}(W)$. 

  Let $\cI_{\Sigma}$ be the subsheaf of 
$\cI$ defined by 
\begin{align*}
\cI_{\Sigma}(U):=\{&b \in \cI(U) \,|\, \text{for any } 
u \in U, \text { there exists } \sigma \in \Sigma \text 
{ such that, for any }\\ &\text{ homomorphism }M_{U, \bar u}
\to \Q^{(+)}_{\ge 0}, \text{ the induced 
} W\times W \to \Q \text{ belongs to }\sigma\}, 
\end{align*}
$U \in (\fs/S)$.  

We define 
$$F'_{\Sig} := F'\times_{\cI}\cI_{\Sigma}$$ 
and prove 
that $F'_{\Sigma}$ is a log smooth log algebraic space in the first sense. 

  To prove it, we may assume that $\Sigma$ is the set of all faces of a cone $\sigma$. 
  In this case, we write $F'_{\Sigma}$ as $F'_{\sigma}$. 
  Let 
$$H_{\sig}:= F'_{\sigma} \times_\cI V_{\sigma} = F'_\cT \times_\cT V_{\sigma}=F'\times_\cI V_{\sigma}.$$ 
  Then $H_{\sig}$ is a log algebraic space in the first sense because 
$F'_\cT$ is so. 

  We consider $F'_{\sigma}$ as a quotient of $H_{\sig}$ by the action by the torus defined as the dense open subspace of $V_{\sig}$ where the log structure is trivial. 

  Let $\oc X$ be the quotient of $\oc H_{\sig}$ by the torus action as a functor from the category of schemes over $\bZ[1/n]$ to the category of sets. 
  We can prove that it is representable by the classical Artin criterion (\cite{Artin:criterion} Theorem 3.4). 
  In fact, the openness of versality is satisfied because $\oc H_{\sig}$ is representable. 
  The effective prorepresentability is also satisfied because $\oc X$ is formally a moduli of log abelian varieties with constant degeneration (\cite{KKN6} Section 2).
  The other conditions are easily verified. 
  Thus we have an algebraic space $\oc X$,  which contains as a dense open subset the part of the moduli of abelian varieties without degeneration. 
  We denote this open set by $U$. 
  Endow $\oc X$ with the log structure by the complement of $U$ and the resulting log algebraic space $X$ in the first sense is log smooth by the local theory.
  Further, the projection $H_{\sig} \to F'_{\sig}$ factors through the natural morphism $H_{\sig} \to X$ because the latter is strict as seen pointwise by the local theory.  
  To see that the induced $X \to F'_{\sig}$ is an isomorphism, by the same argument as in \ref{logreg} with the use of the local theory, it is enough to show that $F'_{\sig}$ has the two properties there. 
  It reduces to that $F'$ has the two properties and again to \ref{1'} and \ref{2'}.  
  Thus $X \to  F'_{\sigma}$ is an isomorphism. 
  We conclude that $F'_{\Sigma}$ is a log smooth log algebraic space in the first sense. 

  By $F'_{\Sigma} = F_{\Sigma} \times_FF'$ and by 
Proposition \ref{sbprop:F'toF}, $F_{\Sigma}$ is also a log smooth log algebraic space in the first sense, which completes the proof of Theorem \ref{thm:mainthm2} except properness. 
\end{para}

\begin{rem}
  The above proof also shows that the space representing $F_{\Sig}$ contains as a dense open subset the moduli space $\cA_{g,n}$ of principally polarized abelian varieties with level $n$ structure. 
\end{rem}

\begin{para}
  {\sc Proof of Theorem \ref{thm:mainthm1} except properness}. 
  We deduce Theorem \ref{thm:mainthm1} except properness from the variant of Theorem \ref{thm:mainthm2} proved in the previous paragraph. 

  We have $F' = \bigcup_{\sig \in K} F'_{\sig}$, where $K$ is as in \ref{cI}. 
  Cover $F'$ by $\coprod_{\sig \in K} F'_{\sig}$, which is a log smooth log algebraic space in the first sense by \ref{F'_sig}.
  For two cones $\tau \subset \sig$, we have that $V_{\tau} \to V_{\sig}$ is log \'etale so that 
$F'_{\tau} \to F'_{\sig}$ is represented by a 
log \'etale morphism of log algebraic spaces in the first sense. 
  This implies that $F'$ is a log smooth log algebraic space in the second sense. 
  Hence, $F$ is also. 
\end{para}

\section{Valuative criterion}
\label{s:val}
  In this section, we prove the properness of our moduli spaces, which completes the proofs of Theorems \ref{thm:mainthm1} and \ref{thm:mainthm2}. 
  We use the valuative criterion and the unique extendability of a log abelian variety over a complete discrete valuation field, which was developed in former parts of our series of papers. 

\begin{para}
\label{properness}
  We use the same notation as before. 
  We have to prove that $F$ and $F_{\Sig}$ are proper over $\bZ[1/n]$. 
  First we prove that $F$ is proper, that is, $F(O_K)=F(K)$ in the notation in \cite{KKN4} 17.3. 
  In order to reduce this equality to the category equivalence $\cC_0^{\mathrm{ptpol}} \simeq \cC_2^{\mathrm{pol}}$ (\cite{KKN6} Theorem 3.4), it is enough to show that for a section $(A,\phi,e) \in F(K)$, where $K$ is a field endowed with an fs log structure, the lattice $Y$ of $A$ and the torus part $T$ of $A$ are unramified and the abelian part $B$ of $A$ is of semistable reduction. 
  To see it, it is enough to show that the local monodromy acts on the Tate module of $A$ trivially, which is by Lemma \ref{lem:aut-triv}. 
  This completes the proof of Theorem \ref{thm:mainthm1}. 
\end{para}

\begin{para}
  Next we prove that $F_{\Sig}$ is proper. 
  The proof is parallel to that of \cite{KKN4} Proposition 11.3, which says that the properness of a weak log abelian variety implies the properness of its model with respect to a complete fan.

  First, to see that $F_{\Sig} \to \Spec(\bZ[1/n])$ is separated, we use the following valuative criterion for separatedness for algebraic spaces. 
\end{para}

\begin{prop}
\label{p:valcri}
  Let $Y$ be a locally noetherian scheme. 
  Let $f\colon X \to Y$ be a morphism of algebraic spaces. 
  Assume 

$(1)$ $f$ is locally of finite type, and 

$(2)$ for any commutative diagram 
$$\begin{CD}
\Spec(K) @>>> X \\
@VVV @VVV \\
\Spec(O_K) @>>> Y, 
\end{CD}
$$
where $K$ is a discrete valuation field with the valuation ring $O_K$, 
there exist at most one morphism 
$\Spec(O_K) \to X$ such that $\Spec(K) \to  X$ coincides 
with $\Spec(K) \to \Spec(O_K) \to X$ and $\Spec(O_K) \to X \to Y$ 
coincides with $\Spec(O_K) \to Y$. 

Then $f$ is separated. 
\end{prop}

  See \cite{SP} Lemma 19.1 (0ARJ) for a proof. 
  Note that, by definition we have adopted, an algebraic space is quasiseparated.

\begin{para}
  Since we already see that $F_{\Sig}$ is log smooth, the condition (1) in Proposition \ref{p:valcri} is satisfied. 
  Let the notation be as in Proposition \ref{p:valcri} (2) with $X=F_{\Sig}$ and $Y=\Spec(\bZ[1/n])$. 
  Endow $\Spec(K)$ with the inverse image log structure from $F_{\Sig}$ and $\Spec(O_K)$ with the direct image log structure from $\Spec(K)$. 
  Then, since a morphism $\Spec(O_K) \to X$ of algebraic spaces making the diagram in Proposition \ref{p:valcri} (2) commute uniquely extends to a morphism of log algebraic spaces making the diagram of log algebraic spaces commute, it is sufficient to show that $F_{\Sig}(O_K) \to F_{\Sig}(K)$ is injective. 
  This is by $F(O_K)=F(K)$ in \ref{properness} and $F_{\Sig} \subset F$. 
  Hence $F_{\Sig}$ is separated. 
\end{para}

\begin{para}
  Now we apply to $F_{\Sig} \to \Spec(\bZ[1/n])$ another valuative criterion \cite{KKN4} Proposition 11.4 for universal closedness.
  Since we already see that $F_{\Sig}$ is separated, it suffices to check the conditions (1) and (2) in \cite{KKN4} Proposition 11.4.

  First $F_{\Sig}$ is quasicompact because, in virtue of the condition (4) in Definition \ref{defn:adm-decomp}, $F_{\Sig}$ is the union of a finite number of $\sigma$-loci ($\sig \in \Sig$), and each $\sigma$-locus is quasicompact by the local theory (\cite{KKN6} Section 2).
  Since we already know that $F_{\Sig}$ is log smooth, it is also of finite type, that is, the condition (1) in \cite{KKN4} Proposition 11.4 is satisfied. 

  The condition (2) in \cite{KKN4} Proposition 11.4 can be deduced from $F(O_K)=F(K)$ in \ref{properness} exactly in the same way as in \cite{KKN4} 11.6--11.7.
  (Note that in the third last line in \cite{KKN4} 11.7, $A^{(\Sig)}(S'') \to A^{(\Sig)}(\eta'')$ should be replaced by $A^{(\Sig)}(\eta'') \to A(\eta'')$.)
  This completes the proof of Theorem \ref{thm:mainthm2}. 
\end{para}

\section{Toroidal compactifications}
\label{s:toroidal}
  Finally we prove Theorem \ref{thm:toroidal}, which says that the toroidal compactification is the underlying space of our space. 

  \begin{para}
{\sc Proof of Theorem \ref{thm:toroidal}}. 
  We give two proofs, which have an outline in common: 
  In both proofs, we construct a family of log abelian varieties over the Faltings--Chai space $\overline {\cA}_{g,n}$, which will coincide with our universal family. 
  The family decides a morphism $i$ from $\overline {\cA}_{g,n}$ to our space because ours is a fine moduli, and we prove that $i$ is an isomorphism. 

  In the first proof, we construct the family of log abelian varieties by gluing the local universal families over our local moduli in \cite{KKN6}.  
  The compatibility of the local families is checked by using the local properties of Faltings--Chai's family $P$ over $\overline {\cA}_{g,n}$, which compactifies the universal family of abelian varieties over the open set ${\cA}_{g,n}$ without degeneration of $\overline {\cA}_{g,n}$. 
  After we prove Theorem \ref{thm:toroidal}, we will show that their family $P$ is a model of our family.

  In the second proof, we construct the family of log abelian varieties directly from $P$ by \lq\lq recovering from models''-procedure. 
  Then we will see that $P$ is a model of our family as soon as the proof of Theorem \ref{thm:toroidal} will finish.

  Notice that these two constructions generalize the two constructions of universal log elliptic curves over $X(N)$ in \cite{KKN3} Section 3 and in \cite{KKN3} Appendix, respectively.
\end{para}

\begin{para}
\label{1stpf-thm:toroidal}
  We explain the first proof. 
  Let $x$ be a point of $\overline {\cA}_{g,n}$. 
  Since the formal completion of $\overline {\cA}_{g,n}$ at $x$ can be identified with our local moduli in \cite{KKN6}, by GAGF for log abelian varieties (\cite{KKN6} Theorem 6.1), we have a family $A_x$ of log abelian varieties (with additional structures, omitted now) over the completion of $\overline {\cA}_{g,n}$ at $x$. 
  We claim that they are compatible with each other and with the universal family of abelian varieties over ${\cA}_{g,n}$.
  Then, they are glued globally and gives a morphism $i$ from $\overline {\cA}_{g,n}$ to our space (cf.\ \cite{KKN2} Proposition 4.8). 

  To prove this compatibility, we use the local portion of the theory of Faltings--Chai.
  Let $P$ be the family constructed by Faltings--Chai associated to an admissible polyhedral cone decomposition $\{\tau_{\beta}\}$.  
  We endow $P$ with the fs log structure determined by the inverse image of $\overline{\cA}_{g,n}-{\cA}_{g,n}$.
  Then, over the formal completion of $\overline {\cA}_{g,n}$ at $x$, which is identified with our local moduli, the pullback of $P$ coincides with a model of our formal universal family associated to the fan decided by $\{\tau_{\beta}\}$. 
  Hence, by GAGF, over the completion at $x$, the pullback of $P$ can be regarded as a model of $A_x$. 
  Then, each $A_x$ is compatible with the universal family of abelian varieties on ${\cA}_{g,n}$ because its model is so.
  To say further that the $A_x$ are compatible with each other because their models are compatible, it is enough to know that, for a complete fan $\Sig'$,  any isomorphism of the $\Sig'$-models of two log abelian varieties uniquely extends to an isomorphism of these log abelian varieties. 
  This is the usual recovering from models and in fact generally valid (see Remark \ref{r:enhance}). 
  Here we use the category equivalence Theorem 4.7 in \cite{KKN6}. 
  To use it, the complete fan $\Sig'$ around $x$ associated with $\{\tau_{\beta}\}$ should be wide. 
  Since this wideness is not necessary satisfied by a general $\{\tau_{\beta}\}$, we work locally around $x$, take, for example, the first standard fan, which is wide (\cite{KKN5} Proposition 10.3), and apply the construction of good compact formal models at the boundary in \cite{FC} Chapter VI 1.11 to this local situation. 
  Then the compatibility of models still holds (see \cite{FC} Chapter VI 1.11 p.206), and by Theorem 4.7 in \cite{KKN6}, we see that the $A_x$ are compatible, as desired. 

  We prove that $i$ is an isomorphism. 
  Since it is formally an isomorphism, by the argument in \cite{KKN3} 5.4, we reduce to the two properties of our functor, that is, finiteness and uniqueness for GAGF. 
Both are valid as our functor is represented by Theorem \ref{thm:mainthm2}. 
  Thus we complete the proof of Theorem \ref{thm:toroidal}.

  Additionally, under the identification of their space and our space, which is just established, $P$ is a model of our universal family because it is so formally around each point, also on 
${\cA}_{g,n}$, and isomorphisms are glued since they are compatible on ${\cA}_{g,n}$.
\end{para}

\begin{para}
\label{2nd}
  Next we explain the second proof.
  By the argument of the second last paragraph of \ref{1stpf-thm:toroidal}, it suffices to construct an appropriate family of log abelian varieties on $\overline {\cA}_{g,n}$. 
  We construct it from $P$ in \ref{1stpf-thm:toroidal} directly.

  On $\overline {\cA}_{g,n}$, there are semiabelian schemes $G$ and $G^t$ and the pairing $\underline X(G) \times \underline X(G^t) \to \Gmlog/\Gm$ 
globally as in \cite{FC} Chapter III, Section 10. Hence the sheaf 
$$Q: = \cHom(\overline X, \Gmlog/\Gm)^{(\overline Y)}/\overline Y$$ is defined globally, where 
$\overline X = \underline X(G)$ and $\overline Y = \underline X(G^t)$.
  The $\{\tau_{\beta}\}$ in \ref{1stpf-thm:toroidal} decides a subsheaf $\Sig'$ of $Q$ coming from complete fans (cf.\ \cite{KKN6} 4.4).

  Then $P$ is a $G$-torsor over $\Sig'$. %
  This can be checked formally and reduces to the fact explained in \ref{1stpf-thm:toroidal} that %
$P$ formally coincides with the model of our universal family. 
  We construct a family of log abelian varieties from $P$. 
  Below, for simplicity, we assume that $\Sig'$ comes from complete and wide fans and use the category equivalence Theorem 4.7 in \cite{KKN6}. 
  The authors do not know if one can take such a $\{\tau_{\beta}\}$.
  But even if it is not the case, we can generalize the above category equivalence to the one without the assumption of the wideness of fans and use it.
  See Remark \ref{r:enhance}.

  Under the above assumption, it is enough to construct an object of the category $\cB$ in \cite{KKN6} Section 4 of the models and it is almost done in \cite{FC} except the construction of the partial group law 
$(P \times P)' \to P$, where $(P \times P)'$ is constructed in \cite{FC} Chapter VI as a compactification of the product of two copies of the universal abelian varieties. 
  To give it, we start with the partial group law in formal situation. 
  By GAGF, it induces a partial group law over the completion of each point of the base. 
  Then, it is glued globally because it is compatible with the group law on the open subset ${\cA}_{g,n}$. 
  Thus we have an object of $\cB$.
\end{para}

\begin{rem}
\label{r:enhance}
  The ``recovering from models''-type statement without wideness of fans, which is mentioned in \ref{1stpf-thm:toroidal} and in \ref{2nd}, is based on a generalization of the k\'et presentation of $\tilde A$ in \cite{KKN6} 4.3 as follows. 
  In there, we consider the morphism $Z:=\coprod_{n \ge0} \tilde A^{(\Sig)} \to \tilde A; (x,n) \mapsto x^{\ell^n}$ (the notation is as in there). 
  Recall that, since it is k\'et surjective by virtue of wideness of $\Sig$, we can recover $\tilde A$ by dividing $Z$ by some relation. 
  If $\Sig$ is not necessarily wide, this surjectivity is no longer valid. 
  However, even in case where $\Sig$ is only complete, the morphism $\coprod_{y \in Y} \coprod_{n \ge0} \tilde A^{(\Sig)} \to \tilde A; (x,n,y) \mapsto x^{\ell^n}y$ is still k\'et surjective. 
  Hence we could formulate the category equivalence also in this case, though the relation becomes more complicated.
  See \cite{KKN2} Section 1 and \cite{KKN3} Appendix for such variants of constructions to recover a log elliptic curve from models, which work without wideness of fans. 
\end{rem}

\noindent Takeshi Kajiwara

\noindent Department of Applied mathematics \\
Faculty of Engineering \\
Yokohama National University \\
Hodogaya-ku, Yokohama 240-8501 \\
Japan

\noindent kajiwara@ynu.ac.jp
\par\bigskip\par

\noindent Kazuya Kato

\noindent 
Department of Mathematics
\\
University of Chicago
\\
5734 S.\ University Avenue
\\
Chicago, Illinois, 60637 \\
USA

\noindent kkato@math.uchicago.edu
\par\bigskip\par

\noindent Chikara Nakayama

\noindent Department of Economics \\ Hitotsubashi University \\
p2-1 Naka, Kunitachi, Tokyo 186-8601 \\ Japan

\noindent c.nakayama@r.hit-u.ac.jp
\end{document}